\documentclass[11pt,letterpaper]{article}
\usepackage[margin=1in]{geometry}
\usepackage[T1]{fontenc}
\usepackage{newtxtext}
\usepackage{microtype,amsmath,amsthm,mathtools}
\usepackage{newtxmath}
\usepackage{bm,mathrsfs}
\DeclareFontFamily{U}{rsfs}{}
\DeclareFontShape{U}{rsfs}{m}{n}{<-6>rsfs5<6-8>rsfs7<8->rsfs10}{}
\usepackage{graphicx,booktabs,array,enumitem,xcolor,fancyhdr}
\usepackage{setspace}
\usepackage[numbers,sort&compress]{natbib}
\usepackage{hyperref}
\hypersetup{colorlinks=true,linkcolor=blue!45!black,citecolor=blue!45!black,urlcolor=blue!45!black}
\renewenvironment{abstract}{\par\noindent\textbf{Abstract.}\ }{\par\medskip}
\setlist{nosep,leftmargin=*}
\numberwithin{equation}{section}
\newtheorem{theorem}{Theorem}[section]
\newtheorem{proposition}[theorem]{Proposition}
\newtheorem{lemma}[theorem]{Lemma}
\newtheorem{corollary}[theorem]{Corollary}
\theoremstyle{definition}
\newtheorem{assumption}[theorem]{Assumption}

\newtheorem{remark}[theorem]{Remark}
\newcommand{\E}{\mathbb E}
\newcommand{\Q}{\mathbb Q}

\newcommand{\I}{\mathcal I}
\newcommand{\A}{\mathcal A}
\newcommand{\M}{\mathcal M}
\newcommand{\D}{\mathcal D}
\newcommand{\C}{\mathcal C}
\newcommand{\Sset}{\mathcal S}
\newcommand{\1}{\mathbf 1}
\newcommand{\dd}{\,\mathrm d}

\title{\LARGE \bf Consumption and Investment with Multiple Job-Switching Opportunities}
\author{Gugyum Ha\thanks{Department of Mathematics, Sogang University, Seoul 04107, Republic of Korea. Email: \href{mailto:ggha@sogang.ac.kr}{ggha@sogang.ac.kr}.}
\and Junkee Jeon\thanks{Department of Applied Mathematics, Kyung Hee University, Yongin-si 17104, Republic of Korea. Email: \href{mailto:junkeejeon@khu.ac.kr}{junkeejeon@khu.ac.kr}.}}
\date{}
\hypersetup{pdftitle={Consumption and Investment with Multiple Job-Switching Opportunities},pdfauthor={Gugyum Ha and Junkee Jeon}}
\begin{document}
\maketitle
\begin{abstract}
We study finite-horizon optimal switching with affine rewards in a consumption--investment model with finitely many occupational regimes. Directed shortest-path costs yield the exact natural borrowing limit. A dual transformation reduces the control problem to a one-dimensional system of variational inequalities with interconnected obstacles. Convexity and a distributional waiting inequality give local spatial $C^{1,1}$ regularity and smooth fit, including destination ties. These estimates establish strong duality and recover optimal consumption, investment, and switching strategies. A reflected backward stochastic differential equation representation yields a set-based integral formula without prescribed boundary geometry. Every connected component of the spatially regular locus is a maximal smooth graph over time. Exact examples exhibit disconnected contact sets, temporal boundaries under strict triangle costs, and a corner in an individual boundary component under additive route costs. A convergent time discretisation completes the analysis.
\end{abstract}

\section{Introduction}\label{sec:introduction}
Occupational decisions affect both the income available for consumption and the risk that an agent can bear through financial investment. A job with a higher wage may require greater effort, while moving between jobs can involve a substantial monetary cost. When several alternatives remain available before retirement, the value of a move depends on the opportunities that follow it. This interaction makes occupational choice a dynamic component of the consumption and investment problem initiated by Merton \cite{Merton69,Merton71}.

We consider an agent who can switch repeatedly among $N$ jobs before a fixed mandatory retirement date. Each job has a constant wage and an additive labour disutility. The agent trades in a complete Black--Scholes market and has common consumption and terminal utility across jobs. Wages and disutilities may coincide or have different orderings, and switching costs may depend on the direction of a move. These features accommodate jobs that cannot be ranked by a single wage--effort index.

\subsection{Main results and approach}
Our first result identifies the financing domain. When finite instantaneous routes are available, directed shortest paths give the effective switching costs. Maximal net human capital is the upper envelope of discounted destination wages less effective entry costs. Theorem~\ref{thm:human} proves that this envelope gives the exact natural borrowing limit for each initial job.

The second result concerns value regularity and control recovery. A martingale dual transformation separates common utility from an $N$-mode switching problem with one continuous state. Convexity, a horizon Lipschitz estimate, and a distributional waiting inequality yield locally bounded second spatial derivatives without excluding destination ties (Theorem~\ref{thm:regular}). Spatial smooth fit therefore holds at every contact. These estimates identify the associated reflected backward equation and establish strong duality, a unique multiplier, and optimal controls throughout the natural wealth domain (Theorem~\ref{thm:duality}).

The third result separates value regularity from boundary geometry. Consolidating active routes permits a continuing destination to be selected and yields a set-based integral representation (Theorem~\ref{thm:integral}). Where the active destination is unique and the gap has a nonzero spatial half-space blow-up, each connected regular component is a maximal $C^\infty$ graph over calendar time (Theorems~\ref{thm:regular_components} and~\ref{thm:cinfty_regular}). Exact examples show why these hypotheses matter: contact sets can be disconnected, strict triangle costs permit temporal boundaries, and additive route costs can produce a corner in an individual boundary component. Spatial smooth fit survives that corner.

\subsection{Related literature}
The financial formulation uses the martingale approach of Karatzas, Lehoczky and Shreve \cite{KLS}, Cox and Huang \cite{CH}, and Karatzas and Shreve \cite{KS}. Cuoco \cite{Cuoco} and Koo \cite{Koo} incorporate labour income and borrowing restrictions; Bodie, Merton and Samuelson \cite{BMS} study flexible labour supply. Discretionary retirement is considered by Karatzas and Wang \cite{KW}, Farhi and Panageas \cite{FP}, Dybvig and Liu \cite{DL}, and Yang and Koo \cite{YK}. Here retirement is fixed, but occupational changes remain available until that date.

The closest occupational-choice models are Jeon and Park \cite{JP}, Shim and Jeon \cite{JS}, and Yang and Jeon \cite{YJ}. The last paper derives a finite-horizon two-job dual problem with a double obstacle and two free boundaries. Ha, Jeon and Ok \cite{HJO} allow time-varying switching costs. We retain constant primitives but allow arbitrarily many jobs, introducing interconnected destinations, asymmetric routes, and simultaneous contacts that cannot generally be analysed as independent pairwise problems.

\paragraph{Finite-horizon switching and reflected BSDEs.}
Tang and Yong \cite{TY} develop a viscosity approach to finite-horizon switching and impulse control. Hamad\`ene and Jeanblanc \cite{HJ} study two-regime starting and stopping, while Djehiche, Hamad\`ene and Popier \cite{DHP} treat finitely many regimes through Snell envelopes and reflected backward stochastic differential equations (RBSDEs). El Asri and Hamad\`ene \cite{EAH} establish Markovian well-posedness for interconnected variational inequalities, and Martyr \cite{Martyr} treats signed switching costs under suitable structural conditions.

The probabilistic connection builds on the scalar RBSDE theory of El Karoui et al.\ \cite{EKKPQ}. With multiple regimes, each value process is reflected above the maximum of the other values net of switching costs. Hu and Tang \cite{HT} develop this oblique-reflection interpretation. Hamad\`ene and Zhang \cite{HZ} allow coupling through generators and obstacles; Chassagneux, Elie and Kharroubi \cite{CEK11} obtain well-posedness without driver monotonicity. Hamad\`ene and Morlais \cite{HM} address time- and state-dependent costs, and Elie and Kharroubi \cite{EK14} give representations for switching that also controls volatility.

Our problem has a common diffusion, affine rewards, and zero terminal switching value. Remark~\ref{rem:rbsde} gives its RBSDE representation. Building on the established switching--RBSDE correspondence, we exploit this structure to obtain the financing domain, spatial estimates, and boundary examples. Positive cycle costs permit route consolidation, but strict triangle inequalities and unique active destinations are not imposed globally.

\paragraph{Value regularity and boundary geometry.}
Pham \cite{Pham} establishes one-dimensional spatial smooth fit. Ly Vath and Pham \cite{VP} obtain explicit two-regime solutions, and Pham, Ly Vath and Zhou \cite{PVZ} study infinite-horizon multiple-regime switching. The terminal date here produces a parabolic system, whose value regularity does not imply smooth contact boundaries.

Our boundary analysis uses the scalar obstacle theory of Friedman \cite{Friedman}, Caffarelli, Petrosyan and Shahgholian \cite{CPS}, Blanchet, Dolbeault and Monneau \cite{BDM}, and Lindgren and Monneau \cite{LM}, together with Kukuljan's higher-order boundary Harnack estimates \cite{Kuk}. Related financial double obstacles arise in Dai and Yi \cite{DY}. These results apply where a single active switching obstacle is smooth; the multiple-destination example shows why that reduction need not extend across all contacts.

\paragraph{Numerical approximation and reinforcement learning.}
The computational approach is related to the switching discretisations of Gassiat, Kharroubi and Pham \cite{GKP} and the discrete-time approximations of obliquely reflected BSDEs by Chassagneux, Elie and Kharroubi \cite{CEK}. Our exact time-grid recursion is monotone and retains all destinations. We prove convergence by approximating admissible switching times, complementing the viscosity-scheme perspective of Barles and Souganidis \cite{BS}.

Huang, Li, Yu and Zhou \cite{HLYZ} study finite-horizon switching through entropy regularisation and reinforcement learning in continuous time. Randomising switching times and destinations leads to an exploratory HJB system, policy-iteration convergence, and convergence to the classical value as the temperature vanishes. Their model-free approach complements our analytical and model-based treatment. However, their bounded rewards, uniform ellipticity, and strict triangle costs differ from our affine rewards, geometric diffusion, and possible triangle equalities. Extending their analysis to this model therefore requires additional work; our explicit examples provide potential benchmarks.

Section~\ref{sec:model} specifies the model and establishes the financing domain. Section~\ref{sec:solution} develops the dual switching system, optimal controls, RBSDE representation, and integral formula. Section~\ref{sec:geometry} studies localisation, smooth regular components, and exact boundary examples. Section~\ref{sec:computation} gives the CRRA formulas, computational results, and the vanishing-cost limit. Section~\ref{sec:conclusion} concludes.

\section{The model and its financing domain}\label{sec:model}
\subsection{Financial market and occupational regimes}
Let $\I=\{1,\ldots,N\}$, where $2\le N<\infty$. On a probability space with the usual augmented Brownian filtration $(\mathcal F_s)_{s\ge0}$, the traded assets satisfy
\[
 \frac{\dd S_s^0}{S_s^0}=r\dd s,\qquad
 \frac{\dd S_s^1}{S_s^1}=\mu\dd s+\sigma\dd B_s,
 \qquad r>0,\quad\sigma>0,\quad\theta=\frac{\mu-r}{\sigma}\ne0.
\]
The retirement date is $T>0$ and the discount rate is $\beta>0$. Job $i$ pays $\epsilon_i\ge0$ and incurs disutility $L_i\ge0$. A common constant can be added to all $L_i$ without changing any switching decision, so the nonnegative convention entails no loss for an arbitrary finite list of real disutilities.

\subsection{Switching costs, preferences, and admissibility}
Let $k_{ij}>0$ for $i\ne j$, and $k_{ii}=0$. A finite route at one decision time is permitted, and all its monetary costs are charged at that time. Define
\begin{equation}\label{eq:distance}
 d_{ii}=0,\qquad
 d_{ij}=\min_{i=i_0,\ldots,i_m=j}\sum_{\ell=1}^m k_{i_{\ell-1}i_\ell},\quad i\ne j.
\end{equation}
Every minimising route can be chosen simple. Consequently, it contains at most $N-1$ edges, $d_{ij}>0$, and
\begin{equation}\label{eq:triangle}
 d_{ij}\le d_{ik}+d_{kj}.
\end{equation}
An effective switch $i\to j$ means executing such a route and paying $d_{ij}$.

\begin{proposition}\label{prop:closure}
The model allowing finite instantaneous routes with costs $k$ and the model with effective direct costs $d$ have the same value functions. Every original route can be replaced by an effective switch with the same endpoints and no greater monetary cost, and every effective switch can be implemented in the original model.
\end{proposition}
\begin{proof}
Replace each route by a shortest route with the same endpoints. The elapsed time and terminal job are unchanged, whereas the monetary cost weakly decreases. The saved amount can be retained in the risk-free asset, so this replacement cannot reduce attainable utility. Conversely, every effective switch can be implemented by a shortest route in the original graph. Positivity removes cycles and ensures finitely many original edges per effective switch. The two comparisons prove equality of the value functions.
\end{proof}

\begin{remark}\label{rem:timing}
This routing convention matters. If the original model forbids simultaneous successive switches and $k_{ij}>k_{i\ell}+k_{\ell j}$, a cheapest instantaneous route may be unattainable in that model. Vanishing delays can approximate its value, but do not by themselves prove attainment. There is no such distinction when the original costs already satisfy the triangle inequality. Throughout the paper, admissible switches are interpreted under this routing convention.
\end{remark}

A strategy $\alpha=((\tau_n,j_n))_{n\ge1}$ consists of stopping times before $T$ and $\mathcal F_{\tau_n}$-measurable destination jobs $j_n\in\I$. After route consolidation, distinct decision times are strictly increasing. The resulting job process is $I_s$, initially $I_{t-}=i$. Strategies have finitely many switches almost surely and satisfy
\[
 \E_t\sum_n H_{t,\tau_n}d_{j_{n-1}j_n}<\infty,
\]
where $H$ is the state-price density defined in \eqref{eq:H}. No switch is made at $T$, because terminal utility is job-independent and switching costs are positive.

With consumption $c_s>0$ and risky investment $\pi_s$ in dollars, wealth satisfies
\begin{equation}\label{eq:wealth}
 \dd X_s=(rX_s+(\mu-r)\pi_s+\epsilon_{I_s}-c_s)\dd s
       +\sigma\pi_s\dd B_s-\sum_n d_{j_{n-1}j_n}\delta_{\tau_n}(\dd s).
\end{equation}
Here $j_0=i$ and the last term denotes monetary jumps. Consumption and investment are progressively measurable, with $\int_t^T(c_s+\pi_s^2)\dd s<\infty$ almost surely. We impose the natural borrowing floor derived in \eqref{eq:floor}, require $X_T>0$, and consider controls for which the expected utility is well defined. The objective is
\begin{equation}\label{eq:primal}
 V_i(t,w)=\sup_{c,\pi,\alpha}\E_{t,w,i}
 \left[\int_t^T e^{-\beta(s-t)}\{U_1(s,c_s)-L_{I_s}\}\dd s
       +e^{-\beta(T-t)}U_2(T,X_T)\right].
\end{equation}

\begin{assumption}\label{ass:utility}
For $a=1,2$, $U_a(s,\cdot)$ is smooth, strictly increasing, has strictly negative second derivative, and satisfies the Inada limits. Its conjugate and inverse marginal utility are
\[
 \widetilde U_a(s,y)=\sup_{c>0}\{U_a(s,c)-yc\},\qquad
 \mathcal I_a(s,y)=(\partial_cU_a(s,\cdot))^{-1}(y).
\]
The expectations and differentiations in \eqref{eq:R} below are finite and locally justified. A convenient sufficient condition is a uniform bound by $C(1+y^p+y^{-p})$, for some $C,p>0$, on the conjugates and their first three $y$-derivatives. Power and logarithmic utilities satisfy this condition. Only these finiteness and differentiation properties, together with strict concavity and the Inada conditions, are used in the duality argument.
\end{assumption}

\subsection{Net human capital and the natural borrowing limit}
Write
\[
 a_q(u)=\frac{1-e^{-qu}}q\quad(q\ne0),\qquad a_0(u)=u,
 \qquad\tau=T-t.
\]
Define the human-capital envelope
\begin{equation}\label{eq:h}
 {\displaystyle h_i(t)=\max_{j\in\I}\{\epsilon_j a_r(T-t)-d_{ij}\}.}
\end{equation}
The initial job is included through $d_{ii}=0$.

\begin{theorem}\label{thm:human}
For every initial job $i$ and every finite switching schedule, the pathwise discounted income less switching costs, with deterministic discount rate $r$, is at most $h_i(t)$. Equality is attainable by an immediate shortest route to a maximiser in \eqref{eq:h}, followed by no further switching. The exact natural wealth domain is
\begin{equation}\label{eq:floor}
 X_s\ge-h_{I_s}(s),\qquad w>-h_i(t).
\end{equation}
The strict inequality is the domain used for strictly positive consumption and terminal wealth.
\end{theorem}
\begin{proof}
Set $t=0$ and let the remaining horizon be $\tau$. For a realised finite job sequence $j_0=i,j_1,\ldots,j_m$ with switching times $s_1\le\cdots\le s_m$, its discounted cash value equals
\begin{equation}\label{eq:cashpoly}
 \frac{\epsilon_i-\epsilon_{j_m}e^{-r\tau}}r
 +\sum_{n=1}^m e^{-rs_n}
   \left(\frac{\epsilon_{j_n}-\epsilon_{j_{n-1}}}r-d_{j_{n-1}j_n}\right).
\end{equation}
This is linear in $z_n=e^{-rs_n}$ on the compact polytope
\[
 1\ge z_1\ge\cdots\ge z_m\ge e^{-r\tau}.
\]
Every vertex places a prefix of the switches at time zero and the remaining switches at $\tau$. Deleting all switches at $\tau$ increases the cash value. A prefix ending in job $j$ therefore gives at most $\epsilon_j a_r(\tau)-d_{ij}$. Maximising proves the upper bound. Switching immediately to a maximiser and remaining there proves attainability. The argument applies to each realised random schedule, so also to every conditional risk-neutral expectation.

Subsection~\ref{sec:budget} shows that future positive consumption and terminal wealth can be financed exactly on $w>-h_i(t)$ and that any admissible financing strategy obeys the floor. At equality there is no positive budget left beyond maximal net labour income. Thus strict positivity requires the open domain in \eqref{eq:floor}.
\end{proof}

\begin{corollary}\label{cor:cashcritical}
If $\epsilon_j>\epsilon_i$ and $rd_{ij}<\epsilon_j-\epsilon_i$, then immediate switching to $j$ gives more net human capital than remaining in $i$ precisely when
\begin{equation}\label{eq:cashcritical}
 T-t>\tau_{ij}^{\rm cash}:=-\frac1r\log\left(1-\frac{rd_{ij}}{\epsilon_j-\epsilon_i}\right).
\end{equation}
If $rd_{ij}\ge\epsilon_j-\epsilon_i$, it never does so at a finite horizon. If $\epsilon_j\le\epsilon_i$, it cannot increase net human capital.
\end{corollary}
\begin{proof}
Solve $(\epsilon_j-\epsilon_i)a_r(T-t)>d_{ij}$.
\end{proof}
These dates compare borrowing capacity. Utility-optimal switching dates also depend on labour disutilities and subsequent occupational opportunities.

\subsection{Static and dynamic budgets}\label{sec:budget}
For $s\ge t$, put
\begin{equation}\label{eq:H}
 H_{t,s}=\exp\{-r(s-t)-\tfrac12\theta^2(s-t)-\theta(B_s-B_t)\}.
\end{equation}
Let $\Q$ be the risk-neutral measure, so $\E_t[H_{t,s}F_s]=\E_t^{\Q}[e^{-r(s-t)}F_s]$ for integrable $F_s$.

\begin{proposition}\label{prop:budget}
An admissible strategy with initial wealth $w$, terminal wealth $b\ge0$, and the floor \eqref{eq:floor} satisfies
\begin{equation}\label{eq:budget}
 \E_t\left[\int_t^T H_{t,s}(c_s-\epsilon_{I_s})\dd s
        +\sum_n H_{t,\tau_n}d_{j_{n-1}j_n}+H_{t,T}b\right]\le w.
\end{equation}
Conversely, given a switching strategy, positive consumption and nonnegative terminal wealth with finite absolute budget terms, equality in \eqref{eq:budget} is sufficient for financing by a portfolio. At any time immediately after a switch, the financing wealth is
\begin{equation}\label{eq:replication}
 X_s=\E_s^{\Q}\left[\int_s^T e^{-r(u-s)}(c_u-\epsilon_{I_u})\dd u
       +\sum_{\tau_n>s}e^{-r(\tau_n-s)}d_{j_{n-1}j_n}
       +e^{-r(T-s)}b\right].
\end{equation}
This wealth satisfies \eqref{eq:floor}. Immediately before a switch its cost is included in the remaining liability.
\end{proposition}
\begin{proof}
Integration by parts applied to $H_{t,s}X_s$ gives the static budget process plus a local martingale. The lower bound $X_s\ge-h_{I_s}(s)\ge-\epsilon_{\max}/r$, together with the nonnegative consumption and switching liabilities, makes the relevant local martingale a supermartingale after localisation. Its negative part is controlled by an integrable multiple of $\sup_{u\le T}H_{t,u}$ and the bounded income stream. Fatou's lemma gives \eqref{eq:budget}.

For the converse, the conditional expectation of the total discounted liability is a uniformly integrable martingale. The Brownian martingale representation theorem supplies a risky-asset integrand \cite{KS}; subtracting the liabilities already paid gives \eqref{eq:wealth}. Theorem~\ref{thm:human}, applied conditionally to the remaining job schedule, bounds its net future labour income by $h_{I_s}(s)$. Formula \eqref{eq:replication} therefore gives $X_s+h_{I_s}(s)\ge0$. Positivity of the future consumption or bequest makes this inequality strict before the horizon. Local square integrability and terminal convergence follow from the martingale representation; stronger moment assumptions in Assumption~\ref{ass:utility} give the usual integrability of the optimal controls constructed below.

Finally, when $w>-h_i(t)$, first choose a human-capital-maximising immediate route. The remaining positive total budget $w+h_i(t)$ can finance positive deterministic consumption and a positive terminal amount. If $w\le-h_i(t)$, \eqref{eq:budget} and the human-capital bound exclude such a strictly positive allocation.
\end{proof}

\section{The dual problem and optimal strategies}\label{sec:solution}
\subsection{Dual formulation and the switching system}
Define
\begin{equation}\label{eq:Y}
 Y_s^{t,y}=y e^{\beta(s-t)}H_{t,s},\qquad
 \dd Y_s=(\beta-r)Y_s\dd s-\theta Y_s\dd B_s,
\end{equation}
and the discounted generator
\begin{equation}\label{eq:generator}
 \A\phi=\tfrac12\theta^2y^2\phi_{yy}+(\beta-r)y\phi_y-\beta\phi.
\end{equation}
The common utility component is
\begin{equation}\label{eq:R}
 R(t,y)=\E\left[\int_t^T e^{-\beta(s-t)}\widetilde U_1(s,Y_s^{t,y})\dd s
       +e^{-\beta(T-t)}\widetilde U_2(T,Y_T^{t,y})\right].
\end{equation}
It solves $R_t+\A R+\widetilde U_1=0$ with terminal value $\widetilde U_2$. Assumption~\ref{ass:utility} gives
\begin{equation}\label{eq:Rderivatives}
 R_y<0,\quad R_{yy}>0,\quad
 \lim_{y\downarrow0}R_y(t,y)=-\infty,\quad
 \lim_{y\uparrow\infty}R_y(t,y)=0.
\end{equation}
For example, differentiation under the expectation gives
\[
 R_y=-\E_t\left[\int_t^T H_{t,s}\mathcal I_1(s,Y_s)\dd s
                 +H_{t,T}\mathcal I_2(T,Y_T)\right].
\]
Strict convexity follows from that of the conjugates; the endpoint limits follow from Inada and monotone or dominated convergence.

Set $f_i(y)=\epsilon_i y-L_i$ and define
\begin{equation}\label{eq:Q}
 Q_i(t,y)=\sup_\alpha\E\left[\int_t^T e^{-\beta(s-t)}f_{I_s}(Y_s^{t,y})\dd s
       -\sum_n e^{-\beta(\tau_n-t)}d_{j_{n-1}j_n}Y_{\tau_n}^{t,y}\right].
\end{equation}
Then the total dual value is
\begin{equation}\label{eq:J}
 J_i(t,y)=R(t,y)+Q_i(t,y).
\end{equation}
The dual switching problem does not contain the consumption utility parameters. They enter the wealth coordinate through $R$.

For a vector $q$, define
\[
 (\M_iq)(t,y)=\max_{j\ne i}\{q_j(t,y)-d_{ij}y\}.
\]
The variational inequalities are
\begin{equation}\label{eq:QVI}
 {\max\{(Q_i)_t+\A Q_i+f_i,\ \M_iQ-Q_i\}=0,
 \qquad Q_i(T,y)=0,\quad i\in\I.}
\end{equation}

For $N\ge3$, the differences $P_{ij}=Q_j-Q_i$ satisfy
\[
 -d_{ji}y\le P_{ij}\le d_{ij}y,\qquad P_{ij}+P_{jk}=P_{ik}.
\]
The presence of a third destination couples these inequalities: either job in a pair may switch to another mode. We therefore retain the full vector system \eqref{eq:QVI}.

\paragraph{A constant-cost representation.}
Let $z=1/y$ and $v_i(t,z)=zQ_i(t,1/z)$. Under $\Q$, the process $Z_s=1/Y_s$ satisfies
\[
 \dd Z_s=(r-\beta)Z_s\dd s+\theta Z_s\dd B_s^{\Q}.
\]
Changing numeraire in \eqref{eq:Q} yields
\begin{equation}\label{eq:constantcost}
 v_i(t,z)=\sup_\alpha\E^{\Q}\left[
  \int_t^T e^{-r(s-t)}(\epsilon_{I_s}-L_{I_s}Z_s)\dd s
  -\sum_n e^{-r(\tau_n-t)}d_{j_{n-1}j_n}\right].
\end{equation}
This is a finite-mode switching problem with constant positive monetary costs, a common geometric diffusion, and affine running rewards. Standard Snell-envelope and viscosity results apply to this representation \cite{DHP,EAH}. The arguments below also establish the properties needed for control recovery directly.

\subsection{Value-function regularity}
Let $L_*=\min_iL_i$ and $L^*=\max_iL_i$.

\begin{lemma}\label{lem:bounds}
The value function $Q_i$ in \eqref{eq:Q} is finite and, for every
$i\in\I$, $t\in[0,T]$, and $y>0$, satisfies
\begin{align}
 y h_i(t)-L^*a_\beta(T-t)
 &\le Q_i(t,y)\le yh_i(t)-L_*a_\beta(T-t),\label{eq:cashbounds}\\
 \lim_{y\downarrow0}Q_i(t,y)&=-L_*a_\beta(T-t),\label{eq:zero}\\
 \lim_{y\uparrow\infty}\frac{Q_i(t,y)}y=h_i(t).\label{eq:infty}
\end{align}
The function $y\mapsto Q_i(t,y)$ is convex, and $|Q_i(t,y)|\le C(1+y)$.
\end{lemma}

\begin{proof}
For a fixed Brownian-adapted policy, write its objective as
\begin{equation}\label{eq:affine}
 yA_i^\alpha(t)-D_i^\alpha(t),\qquad
 A_i^\alpha=\E_t\left[\int_t^T H_{t,s}\epsilon_{I_s}\dd s
               -\sum_nH_{t,\tau_n}d_{j_{n-1}j_n}\right],
\end{equation}
where $D_i^\alpha=\E_t\int_t^T e^{-\beta(s-t)}L_{I_s}\dd s$.
Bounded wages and integrable discounted switching costs make
$A_i^\alpha(t)$ finite. The change-of-measure identity in
Subsection~\ref{sec:budget} and the pathwise bound in
Theorem~\ref{thm:human} give
$A_i^\alpha(t)\le h_i(t)$ for every admissible policy $\alpha$,
with equality attained by an immediate shortest route to a maximiser
in \eqref{eq:h}, followed by no further switching. Moreover, since
$L_*\le L_{I_s}\le L^*$,
\[
    L_*a_\beta(T-t)
    \le D_i^\alpha(t)
    \le L^*a_\beta(T-t).
\]
The upper bounds on $A_i^\alpha$ and $-D_i^\alpha$ give the upper
inequality in \eqref{eq:cashbounds}. Evaluating the payoff at the
cash-maximising policy gives the lower inequality. Thus $Q_i(t,y)$
is finite.

Since the admissible policy set is independent of $y$ and each
fixed-policy payoff in \eqref{eq:affine} is affine in $y$,
$Q_i(t,\cdot)$ is convex.

Dividing \eqref{eq:cashbounds} by $y$ and letting $y\uparrow\infty$
gives \eqref{eq:infty}.
To prove \eqref{eq:zero}, switch immediately to any job with disutility
$L_*$ for the lower bound, and use the upper bound in
\eqref{eq:cashbounds}.
For $t<T$, choosing $j_*$ with $L_{j_*}=L_*$ gives explicitly
\[
 Q_i(t,y)\ge
 y\{\epsilon_{j_*}a_r(T-t)-d_{ij_*}\}-L_*a_\beta(T-t).
\]
Together with \eqref{eq:cashbounds}, this proves \eqref{eq:zero};
at $t=T$ the conclusion follows from the zero terminal value.
Finally, the linear-growth bound follows from \eqref{eq:cashbounds}.
\end{proof}

For an integer $m\ge0$, let $Q_i^{(m)}$ denote the value function obtained by
restricting the admissible strategies in \eqref{eq:Q} to those with
at most $m$ switches, and set
\[
    q_i^{(m)}(\tau,y):=Q_i^{(m)}(T-\tau,y).
\]

\begin{lemma}\label{lem:time}
For every integer $m\ge0$, $y>0$, and $0\le t\le t+\delta\le T$,
\begin{equation}\label{eq:timelip}
    |Q_i^{(m)}(t+\delta,y)-Q_i^{(m)}(t,y)|
    \le C(1+y)\delta,
\end{equation}
where $C$ is independent of $m$. Moreover,
\[
    |Q_i^{(m)}(t,y)|\le C(1+y),
\]
with $C$ independent of $m$. The function
$Q_i^{(m)}(t,\cdot)$ is locally Lipschitz on $(0,\infty)$, and
\begin{equation}\label{eq:gradbound}
    |\partial_y Q_i^{(m)}(t,y)|
    \le C(1+y^{-1})
\end{equation}
for almost every $y>0$, with $C$ independent of $m$.
\end{lemma}

\begin{proof}
Work in elapsed time on a common Brownian space, starting the state
process at $\widehat Y_0=y$. Time homogeneity identifies
$q_i^{(m)}(\tau,y)$ with the value over a horizon of length $\tau$.
For $0\le\delta\le\tau$, truncate any policy for this horizon by
deleting every switch at or after $\tau-\delta$. The resulting policy
is admissible for the shorter horizon, including when a switching
time equals its terminal date. The two policies agree before that
date. Since all deleted switching payments are nonpositive rewards,
the original payoff minus the truncated payoff is at most
\[
 \E\int_{\tau-\delta}^{\tau}e^{-\beta s}
       \max_{k\in\I}|f_k(\widehat Y_s)|\dd s
 \le C(1+y)\delta,
\]
because $f_k(y)=\epsilon_k y-L_k$ and
$\E[\widehat Y_s]=ye^{(\beta-r)s}$.
Conversely, extend a shorter-horizon policy without further switching.
The added running payoff has the same absolute bound. Taking
suprema gives \eqref{eq:timelip}, uniformly in $m$. Both operations
preserve adaptedness on the common elapsed-time Brownian space.

The no-switch policy and Lemma~\ref{lem:bounds} give
\[
 \epsilon_i y a_r(T-t)-L_i a_\beta(T-t)
 \le Q_i^{(m)}(t,y)\le Q_i(t,y)\le C(1+y),
\]
which proves the uniform linear-growth bound.

For each fixed $t$, the admissible strategy set defining
$Q_i^{(m)}$ is independent of $y$, and the payoff under each fixed
strategy is affine in $y$. Hence $Q_i^{(m)}(t,\cdot)$ is finite and
convex on $(0,\infty)$, so it is locally Lipschitz and differentiable
almost everywhere. At every point of differentiability,
\[
    \frac{Q_i^{(m)}(t,y)-Q_i^{(m)}(t,y/2)}{y/2}
    \le \partial_y Q_i^{(m)}(t,y)
    \le
    \frac{Q_i^{(m)}(t,2y)-Q_i^{(m)}(t,y)}{y}.
\]
The uniform linear-growth bound then gives
\eqref{eq:gradbound}, with $C$ independent of $m$.
\end{proof}

\begin{theorem}\label{thm:regular}
The vector $Q$ is the unique continuous viscosity solution of
\eqref{eq:QVI} in the linear-growth class. On compact subsets of
$[0,T)\times(0,\infty)$,
\begin{equation}\label{eq:regular}
 Q_i\in W^{1,2}_{\infty,\mathrm{loc}},\qquad
 Q_i\in C^{0,1},\qquad
 0\le (Q_i)_{yy}\le C\frac{1+y}{y^2}\quad\text{a.e.}
\end{equation}
Moreover,
\[
    Q_i(t,\cdot)\in C^{1,1}_{\mathrm{loc}}(0,\infty),
    \qquad t\in[0,T),
\]
and $(Q_i)_y$ is locally $1/2$-H\"older continuous in $t$.
The spatial Lipschitz and temporal H\"older constants are uniform on
compact cylinders; in particular, $(Q_i)_y$ is jointly continuous.
In its continuation set, $Q_i$ is smooth. No threshold structure,
strict triangle inequality, or absence of destination ties is assumed.
\end{theorem}

\begin{proof}
\emph{Continuity and viscosity characterisation.}
For $m=0$,
\[
 Q_i^{(0)}(t,y)=\epsilon_i a_r(T-t)y-L_i a_\beta(T-t).
\]
Successive optimal-stopping problems give
\begin{align}\label{eq:snellm}
 Q_i^{(m+1)}(t,y)=\sup_{\rho\in[t,T]}\E\Big[&
 \int_t^\rho e^{-\beta(s-t)}f_i(Y_s)\dd s\\
 &+\1_{\{\rho<T\}}e^{-\beta(\rho-t)}
       \max_{j\ne i}\{Q_j^{(m)}(\rho,Y_\rho)-d_{ij}Y_\rho\}\Big].
       \nonumber
\end{align}
The zero terminal reward dominates the payment $-d_{ij}y$ from a
terminal switch. Lemma~\ref{lem:time} and the moments of $Y$ ensure
integrability of the running and stopping rewards, so these Snell
envelopes are well defined. Enlarging the policy class gives
$Q_i^{(m)}\le Q_i^{(m+1)}\le Q_i$.
For the reverse inequality, fix an admissible strategy
$\alpha_0=((\tau_n,j_n))_{n\ge1}$ and define $\alpha_0^{(m)}$ by retaining
the first $m$ switches of $\alpha_0$ and making no further switch
afterward. Since $\alpha_0$ has finitely many switches almost surely,
$\alpha_0^{(m)}$ coincides with $\alpha_0$ on $[t,T]$ for all sufficiently
large $m$ on almost every sample path. Moreover, the payoff random variables corresponding to
$\alpha_0^{(m)}$ are dominated in absolute value by
\[
    C\int_t^T(1+Y_s^{t,y})\dd s
    +
    \sum_n e^{-\beta(\tau_n-t)}
    d_{j_{n-1}j_n}Y_{\tau_n}^{t,y}.
\]
The first term is integrable since
$\E[Y_s^{t,y}]=y e^{(\beta-r)(s-t)}$, while the second is integrable
by the discounted-cost integrability required in the admissibility of $\alpha_0$.
Dominated convergence gives convergence of the truncated payoffs.
Since $\alpha_0^{(m)}$ is admissible for $Q_i^{(m)}$, the payoff of
$\alpha_0$ is at most $\lim_m Q_i^{(m)}(t,y)$. Taking the supremum
over $\alpha_0$ proves $Q_i^{(m)}(t,y)\uparrow Q_i(t,y)$.
The $m$-independent bounds in Lemma~\ref{lem:time} give local
boundedness and equicontinuity of $\{Q_i^{(m)}\}$. Hence
Arzel\`a--Ascoli, together with the pointwise convergence above,
gives locally uniform convergence to $Q_i$ and the continuity of
$Q_i$. The same uniform time and spatial Lipschitz estimates pass to
the limit $Q_i$.

For each $m$, \eqref{eq:snellm} is an optimal-stopping problem, so its
viscosity equation is obtained from \eqref{eq:QVI} with
$Q_i$ replaced by $Q_i^{(m+1)}$ and $\M_iQ$ by
$\M_iQ^{(m)}$. Since both $Q^{(m+1)}$ and $Q^{(m)}$ converge locally
uniformly to $Q$, viscosity stability \cite{CIL} shows that $Q$
satisfies \eqref{eq:QVI}.
For uniqueness, apply the same change of variables as in
\eqref{eq:constantcost} to any candidate solution:
$v_i(t,z)=zQ_i(t,1/z)$. The linear-growth bound is preserved, since
$z(1+1/z)=1+z$. The transformed generator, including discounting, is
\[
 \widehat\A\phi
 =\tfrac12\theta^2z^2\phi_{zz}+(r-\beta)z\phi_z-r\phi.
\]
The switching costs are constant and bounded below by
$d_{\min}:=\min_{i\ne j}d_{ij}>0$, and the running rewards are
$\epsilon_i-L_i z$. These are the structural conditions used in the
comparison argument of \cite{EAH}.

To account for the state domain $(0,\infty)$, localise comparison
to compact subintervals using the common barrier
\[
 \Psi(t,z)=e^{\lambda(T-t)}(1+z^2+z^{-1}).
\]
For sufficiently large $\lambda$,
$-\Psi_t-\widehat\A\Psi>0$. Adding $\eta\Psi$, $\eta>0$, to
every component of a supersolution preserves its obstacle
inequalities and makes the differential inequality strict.
The barrier dominates any difference of linear-growth candidates
as $z\downarrow0$ and $z\uparrow\infty$. Thus comparison on
expanding compact subintervals, by the argument of \cite{EAH},
requires no boundary data at either inaccessible endpoint.
Letting $\eta\downarrow0$ proves uniqueness on $(0,\infty)$.

\emph{A distributional waiting inequality.}
It is always feasible to remain in job $i$ for a deterministic duration
$h$ and then optimise. Thus
\begin{equation}\label{eq:wait}
 Q_i(t,y)\ge\E\left[
 \int_t^{t+h}e^{-\beta(s-t)}f_i(Y_s)\dd s
 +e^{-\beta h}Q_i(t+h,Y_{t+h})\right].
\end{equation}
Multiply by a nonnegative smooth space-time test function compactly
supported away from $y=0$, integrate, transfer the diffusion semigroup
to its adjoint, divide by $h$, and let $h\downarrow0$. Continuity and
local bounds justify the limit and give
\begin{equation}\label{eq:distribution}
 -(Q_i)_t-\tfrac12\theta^2y^2(Q_i)_{yy}
 -(\beta-r)y(Q_i)_y+\beta Q_i-f_i\ge0
\end{equation}
in distributions.

\emph{No singular spatial curvature.}
The uniform estimates from Lemma~\ref{lem:time}, passed to the limit
above, show that $Q_i$ is locally Lipschitz in both $t$ and $y$.
Hence its weak first derivatives satisfy
\[
    (Q_i)_t,\ (Q_i)_y\in L^\infty_{\rm loc}.
\]
Convexity also makes $(Q_i)_{yy}$ a nonnegative distribution, hence a
nonnegative measure.
Rearranging \eqref{eq:distribution} gives the measure inequality
\[
    0\le
    \tfrac12\theta^2y^2(Q_i)_{yy}
    \le
    \bigl[-(Q_i)_t-(\beta-r)y(Q_i)_y+\beta Q_i-f_i\bigr]
    \,\mathrm dt\,\mathrm dy.
\]
The right-hand side has a locally bounded density with respect to
space-time Lebesgue measure. Since $\theta^2y^2$ is bounded away
from zero on compact state intervals, $(Q_i)_{yy}$ is absolutely
continuous and has no singular part there.
Moreover, the time and spatial derivative estimates inherited from
Lemma~\ref{lem:time}, together with
$|Q_i(t,y)|\le C(1+y)$ and the linear growth of $f_i$, yield
\[
    0\le (Q_i)_{yy}
    \le C\frac{1+y}{y^2}
    \quad\text{a.e.}
\]
Together with the bounds for $(Q_i)_t$ and $(Q_i)_y$, this proves
$Q_i\in W^{1,2}_{\infty,\mathrm{loc}}$ and the remaining assertions in
\eqref{eq:regular}.

On every smaller compact cylinder, Fubini's theorem and the weak
derivative bound give, for almost every $t$ and every $y,h$ whose
spatial segment stays in the cylinder,
\[
    0\le
    Q_i(t,y+h)+Q_i(t,y-h)-2Q_i(t,y)
    \le Ch^2.
\]
The constant is uniform in $t$. Approximating an arbitrary time by
times in this full-measure set and using continuity of $Q_i$ extends
the inequality to every time slice, including $t=0$.
Each spatial slice is consequently $C^{1,1}_{\rm loc}$.
On a smaller compact cylinder, the spatial $C^{1,1}$ bound controls
the error between $(Q_i)_y$ and a spatial difference quotient by
$Ch$, while the time-Lipschitz estimate controls the difference of
these quotients by $C|t-s|/h$. Hence
\[
 |(Q_i)_y(t,y)-(Q_i)_y(s,y)|
 \le C\{h+|t-s|/h\}.
\]
Taking $h=|t-s|^{1/2}$ for sufficiently close times proves the
time modulus. Together with the uniform spatial Lipschitz bound
for $(Q_i)_y$, this proves joint continuity.
Since $Q$ is a viscosity solution, on the open continuation set, the
strict obstacle inequality reduces \eqref{eq:QVI} to its linear
equation.
This linear equation has smooth coefficients
and smooth source, so interior parabolic regularity gives smoothness
\cite{Lieberman}.
\end{proof}
\subsection{Contact sets and optimal switching}
Define
\begin{align}
 \C_i&=\{(t,y):Q_i(t,y)>\M_iQ(t,y)\},\label{eq:C}\\
 \Sset_i&=\{(t,y):Q_i(t,y)=\M_iQ(t,y)\},\label{eq:S}\\
 \Sset_{ij}&=\{(t,y):Q_i(t,y)=Q_j(t,y)-d_{ij}y\},\quad j\ne i.\label{eq:Sij}
\end{align}
The sets $\C_i$ are open, the contact sets are relatively closed, and $\Sset_i=\bigcup_{j\ne i}\Sset_{ij}$. Pairwise sets may overlap when destinations tie.

\begin{lemma}\label{lem:terminaldest}
At every point of $\Sset_i$, there is a destination $j$ such that
\begin{equation}\label{eq:terminaldest}
 Q_i=Q_j-d_{ij}y,\qquad (t,y)\in\C_j.
\end{equation}
A measurable choice is obtained by selecting the smallest index among such destinations. If all triangles are strict for distinct indices, every active direct destination is already in its continuation set.
\end{lemma}
\begin{proof}
Follow active contact edges. A cycle would imply, after adding its value-matching equalities, $0=y\sum d_{ab}>0$. Hence a chain must end in a continuation mode within at most $N-1$ edges. If its accumulated cost is $K$, then $Q_i=Q_j-Ky$ and $d_{ij}\le K$. The obstacle inequality $Q_i\ge Q_j-d_{ij}y$ forces $K\le d_{ij}$, so $K=d_{ij}$. A finite comparison of continuous functions gives a Borel selector. Under a strict triangle inequality, two consecutive active edges contradict the direct obstacle inequality.
\end{proof}

\begin{theorem}\label{thm:switch}
Starting from $(t,y,i)$, remain in the current job until the first hit of its contact set. At that time choose a continuation destination using Lemma~\ref{lem:terminaldest}, execute a shortest route, and repeat. This policy, denoted $\alpha^{t,y,i,*}$, has finitely many switches almost surely, has finite expected discounted cost, and attains \eqref{eq:Q}.
\end{theorem}
\begin{proof}
After every effective switch the new job has a strictly positive continuation gap. Continuity of $Y$ and $Q$ therefore gives a strictly positive pathwise waiting time until the next switch. Suppose nevertheless that infinitely many switch times accumulate at $s_*\le T$. The path of $Y$ on $[t,T]$ is contained in a compact subset of $(0,\infty)$. Among each block of $N$ transitions there is a cycle, and a subsequence repeats one finite cycle. Passing its contact equalities to the common limit $(s_*,Y_{s_*})$ gives $0=Y_{s_*}\sum d_{ab}>0$, a contradiction. This argument also excludes accumulation at $T$; there $Q_i(T,y)=0$ and positive costs make contact impossible.

Apply the localised It\^o--Krylov formula \cite{Krylov} on the successive continuation intervals. The continuation PDE makes
\[
 e^{-\beta(s-t)}Q_{I_s}(s,Y_s)
 +\int_t^s e^{-\beta(u-t)}f_{I_u}(Y_u)\dd u
 -\sum_{\tau_n\le s}e^{-\beta(\tau_n-t)}d_{j_{n-1}j_n}Y_{\tau_n}
\]
a local martingale under the candidate policy. At a switch, value matching cancels the jump. For any other policy it is a local supermartingale because both QVI inequalities have the appropriate sign. Stop after $m$ switches and in compact state intervals. The bound $|Q_i|\le C(1+Y)$ and the integrable supremum of the geometric diffusion control the value and running-reward terms. Nonnegative accumulated costs then have a uniformly bounded expectation. Monotone convergence for costs and dominated convergence for the other terms remove both localisations. Since only finitely many switches occur on each path and the terminal switching value is zero, equality holds for the candidate and the upper bound holds for every admissible competitor.
\end{proof}

\begin{corollary}\label{cor:smoothfit}
At every contact point in $\Sset_{ij}$,
\begin{equation}\label{eq:smoothfit}
 (Q_i)_y=(Q_j)_y-d_{ij}.
\end{equation}
Thus spatial smooth fit holds even at points with several active destinations. This does not imply that the contact-set boundaries themselves are smooth curves.
\end{corollary}
\begin{proof}
For fixed $t$, the differentiable function $Q_i-Q_j+d_{ij}y$ is nonnegative and vanishes at the contact point. Its derivative is therefore zero.
\end{proof}

\subsection{Strong duality and optimal consumption and investment}
\begin{lemma}\label{lem:envelope}
For any optimal policy in Theorem~\ref{thm:switch}, its cash slope in \eqref{eq:affine} equals
\begin{equation}\label{eq:envelope}
 (Q_i)_y(t,y)=A_i^{\alpha^{t,y,i,*}}(t).
\end{equation}
Also,
\begin{equation}\label{eq:derivlimit}
 \lim_{y\uparrow\infty}(Q_i)_y(t,y)=h_i(t).
\end{equation}
\end{lemma}
\begin{proof}
Hold the optimal Brownian-adapted policy fixed while varying the initial dual scale to $z$. Then $Q_i(t,z)\ge zA^*-D^*$, with equality at $z=y$. Differentiability proves \eqref{eq:envelope}; no differentiation of the switching times is required. For a differentiable convex function, a finite asymptotic slope $Q_i(t,y)/y\to h_i(t)$ implies convergence of its derivative to the same limit. This proves \eqref{eq:derivlimit}.
\end{proof}

\begin{theorem}\label{thm:duality}
Under Assumption~\ref{ass:utility}, for every $w>-h_i(t)$,
\begin{equation}\label{eq:duality}
 {V_i(t,w)=\min_{y>0}\{R(t,y)+Q_i(t,y)+wy\}.}
\end{equation}
The minimiser $y_i(t,w)$ is unique and is characterised by
\begin{equation}\label{eq:wealthmap}
 w=\mathcal W_i(t,y_i(t,w)),\qquad
 \mathcal W_i(t,y):=-R_y(t,y)-(Q_i)_y(t,y).
\end{equation}
For each $t,i$, $\mathcal W_i(t,\cdot)$ is a continuous strictly decreasing bijection from $(0,\infty)$ onto $(-h_i(t),\infty)$.

Let $Y^*$ start at $y_i(t,w)$, and use $\alpha^{t,y_i(t,w),i,*}$ for the job process $I^*$. Optimal consumption and terminal wealth are
\begin{equation}\label{eq:cstar}
 c_s^*=\mathcal I_1(s,Y_s^*),\qquad X_T^*=\mathcal I_2(T,Y_T^*).
\end{equation}
At all post-switch times, optimal wealth is
\begin{equation}\label{eq:xstar}
 X_s^*=-\partial_yJ_{I_s^*}(s,Y_s^*).
\end{equation}
An optimal risky portfolio is, up to $\dd s\otimes\dd\mathbb P$ equivalence,
\begin{equation}\label{eq:pistar}
 {\displaystyle
 \pi_s^*=\frac{\theta}{\sigma}Y_s^*\,
       \partial_{yy}J_{I_s^*}(s,Y_s^*).}
\end{equation}
The second derivative is taken in the current continuation regime; values on switching instants do not affect the portfolio integral.
\end{theorem}
\begin{proof}
\emph{Weak duality.} Apply the conjugate inequality $U_a(c)\le\widetilde U_a(y)+yc$ at $Y_s^{t,y}$ and combine it with \eqref{eq:budget}. Taking the supremum over primal controls gives $V_i(t,w)\le R(t,y)+Q_i(t,y)+wy$ for every $y>0$.

\emph{A unique multiplier.} The function $R+Q_i$ is strictly convex, since $R_{yy}>0$ and $Q_i$ is convex. Its first derivative is continuous. By \eqref{eq:Rderivatives}, \eqref{eq:derivlimit}, and the convex upper derivative bound $(Q_i)_y\le h_i(t)$, its derivative tends to $-\infty$ at zero and to $h_i(t)$ at infinity. Hence $\partial_yJ_i+w$ has exactly one zero whenever $w>-h_i(t)$. This proves the bijection and uniqueness.

\emph{Budget equality and attainment.} Use the optimal switching policy at this multiplier. Equations \eqref{eq:Rderivatives} and \eqref{eq:envelope} imply
\begin{align*}
 -\partial_yJ_i(t,y)=\E_t\Big[&\int_t^T H_{t,s}
      \{\mathcal I_1(s,Y_s)-\epsilon_{I_s^*}\}\dd s\\
     &+\sum_nH_{t,\tau_n}d_{j_{n-1}j_n}
      +H_{t,T}\mathcal I_2(T,Y_T)\Big]=w.
\end{align*}
Proposition~\ref{prop:budget} provides an admissible financing portfolio. Both conjugate inequalities are equalities under \eqref{eq:cstar}, and the switching objective is attained, so weak duality becomes equality. Applying the same conditional argument after any time gives \eqref{eq:xstar}.

\emph{Portfolio and jumps.} Between switches the current dual value is smooth. The diffusion coefficient of $-\partial_yJ_i(s,Y_s^*)$ is $\theta Y_s^*\partial_{yy}J_i(s,Y_s^*)$. Comparison with \eqref{eq:wealth} gives \eqref{eq:pistar}. The almost-everywhere formula can also be obtained from the martingale representation; the local derivative bounds and utility moment bounds justify localisation. At a switch $i\to j$, \eqref{eq:smoothfit} yields
\begin{equation}\label{eq:wealthjump}
 \mathcal W_j(s,Y_s^*)=\mathcal W_i(s,Y_s^*)-d_{ij},
\end{equation}
which is exactly the wealth jump in \eqref{eq:wealth}. The dual state remains continuous. Thus the constructed controls are feasible throughout, and attain the primal value.
\end{proof}

\begin{corollary}\label{cor:wealthcontact}
The wealth switching set from $i$ to $j$ is the image
\[
 \Sset_{ij}^{w}=\{(t,\mathcal W_i(t,y)):(t,y)\in\Sset_{ij}\}.
\]
At such a point,
\[
 V_i(t,w)=V_j(t,w-d_{ij}),\qquad
 (V_i)_w(t,w)=(V_j)_w(t,w-d_{ij})=y.
\]
Since consumption utility is common to all jobs, optimal consumption has no instantaneous jump at an effective job switch.
\end{corollary}
\begin{proof}
Value matching in dual coordinates and \eqref{eq:wealthjump} give a common minimiser in the two Legendre transforms. The consumption statement follows from \eqref{eq:cstar} and continuity of $Y^*$.
\end{proof}

In a smooth primal continuation region, the equivalent HJB system is
\begin{equation}\label{eq:primalHJB}
 \max\left\{
 \begin{aligned}
 &(V_i)_t+(rw+\epsilon_i)(V_i)_w-\beta V_i-L_i
   +\widetilde U_1(t,(V_i)_w)
   -\frac{\theta^2}{2}\frac{((V_i)_w)^2}{(V_i)_{ww}},\\
 &\max_{j\ne i}\{V_j(t,w-d_{ij})-V_i(t,w)\}
 \end{aligned}\right\}=0.
\end{equation}
Only feasible target wealth states are included in the second maximum. The dual formulation avoids solving this nonlinear system directly.

\subsection{Integral representation of the switching regions}
Early-exercise-premium representations connect optimal stopping values to occupation of the stopping region \cite{Peskir,PS}. Here an analogous representation uses a continuing destination at each switching point. Choose this destination according to Lemma~\ref{lem:terminaldest} at every point of $\Sset_i$. Let $\D_{ij}\subseteq\Sset_{ij}$ be its Borel selection sets. For each $i$ these sets are disjoint, cover $\Sset_i$, and satisfy
\[
 (t,y)\in\D_{ij}\quad\Longrightarrow\quad(t,y)\in\C_j.
\]
Define the local switching gain
\begin{equation}\label{eq:g}
 g_{ij}(y)=(\epsilon_j-\epsilon_i-rd_{ij})y-(L_j-L_i).
\end{equation}

\begin{lemma}\label{lem:source}
Almost everywhere,
\begin{equation}\label{eq:source}
 -\{(Q_i)_t+\A Q_i+f_i\}
       =\sum_{j\ne i}g_{ij}(y)\1_{\D_{ij}}(t,y),
 \qquad g_{ij}\ge0\ \text{a.e. on }\D_{ij}.
\end{equation}
\end{lemma}
\begin{proof}
The left side is zero in $\C_i$. On $\D_{ij}$, the function $Q_i-Q_j+d_{ij}y$ is zero and belongs locally to $W^{1,2}_p$ for every finite $p$. The level-set property of weak derivatives gives equality of its time derivatives and its first and second spatial derivatives to zero almost everywhere on that set. Since $j$ is in continuation, $(Q_j)_t+\A Q_j=-f_j$. Also $\A(d_{ij}y)=-rd_{ij}y$. Therefore $(Q_i)_t+\A Q_i=-f_j+rd_{ij}y$, which proves \eqref{eq:source}. The QVI gives its nonnegative sign. The argument is valid for measurable contact sets and does not discard a presumed smooth boundary.
\end{proof}

\begin{remark}[Reflected BSDE representation]\label{rem:rbsde}
The value functions give an explicit instance of the switching--RBSDE
correspondence in \cite{DHP,HT,HZ,HM}. Fix $(t,y)$, write
$Y_s=Y_s^{t,y}$, and define
\[
 \mathscr Y_s^i=Q_i(s,Y_s),\qquad
 \mathscr Z_s^i=-\theta Y_s(Q_i)_y(s,Y_s),\qquad
 K_s^i=\int_t^s\mathfrak r_i(u,Y_u)\dd u,
\]
where
\[
 \mathfrak r_i(u,z)
 =\sum_{j\ne i}g_{ij}(z)\1_{\D_{ij}}(u,z)
 =-\{(Q_i)_t+\A Q_i+f_i\}(u,z)
 \quad\text{a.e.}
\]
A nonnegative Borel representative may be chosen for
$\mathfrak r_i$. Changing it on a space-time null set does not change
$K^i$, since $Y$ has a transition density at every positive elapsed
time. For $t\le s\le T$, the processes satisfy
\begin{align}
 \mathscr Y_s^i
 &=\int_s^T\{f_i(Y_u)-\beta\mathscr Y_u^i\}\dd u
   +K_T^i-K_s^i-\int_s^T\mathscr Z_u^i\dd B_u,
   \label{eq:rbsde}\\
 \mathscr Y_s^i
 &\ge\max_{j\ne i}\{\mathscr Y_s^j-d_{ij}Y_s\},
 \label{eq:rbsde_obstacle}\\
 0&=\int_t^T\left(\mathscr Y_s^i-
   \max_{j\ne i}\{\mathscr Y_s^j-d_{ij}Y_s\}\right)\dd K_s^i
 .\label{eq:rbsde_minimality}
\end{align}
Here $\mathscr Y_T^i=0$, $K_t^i=0$, and each $K^i$ is continuous
and nondecreasing. The localised
It\^o--Krylov formula and Lemma~\ref{lem:source} give
\eqref{eq:rbsde}; the value and derivative estimates, together with
the moments of $Y$, justify passage to $T$ and give square-integrable
value, martingale, and reflection processes. The obstacle inequality
is exactly \eqref{eq:QVI}, and $\mathfrak r_i=0$ in $\C_i$ gives
the Skorokhod minimality condition \eqref{eq:rbsde_minimality}.

The reflection process enforces the interconnected value constraints;
it is not the accumulated switching payment, which has jumps along
a chosen strategy. Switching occurs when a value process meets its
obstacle, with active routes consolidated as in
Lemma~\ref{lem:terminaldest}. Taking discounted expectations in
\eqref{eq:rbsde} yields the integral representation below. Thus the
RBSDE, the QVI, and the switching-premium formula describe the same
value vector, without requiring strict triangle inequalities.
\end{remark}

\begin{theorem}\label{thm:integral}
For every $i$,
\begin{equation}\label{eq:EEP}
 {\begin{aligned}
 Q_i(t,y)={}&\epsilon_i a_r(T-t)y-L_i a_\beta(T-t)\\
 &+\sum_{j\ne i}\int_t^T e^{-\beta(s-t)}
       \E\left[g_{ij}(Y_s^{t,y})\1_{\D_{ij}(s)}(Y_s^{t,y})\right]\dd s.
 \end{aligned}}
\end{equation}
Here $\D_{ij}(s)$ denotes the spatial section of the selected switching set. Together with
\begin{align}
 Q_i&\ge Q_j-d_{ij}y\quad(j\ne i),\label{eq:comp1}\\
 \D_{ij}&\subseteq\{Q_i=Q_j-d_{ij}y\}\cap\C_j,\qquad
 \C_i=(\bigcup_{j\ne i}\D_{ij})^c,\label{eq:comp2}\\
 Q_i&>\M_iQ\ \text{on }\C_i,\qquad
 g_{ij}\ge0\ \text{a.e. on }\D_{ij},\label{eq:comp3}
\end{align}
this is an exact set-based nonlinear Volterra characterisation. The resulting value vector is unique; the selector need not be unique on destination ties.
\end{theorem}
\begin{proof}
Equation \eqref{eq:source} is a linear parabolic equation with a measurable source of at most linear growth and terminal value zero. The localised Feynman--Kac formula, followed by the moment bounds for $Y$, yields \eqref{eq:EEP}. The baseline term is the expectation of $f_i$ when no premium is present. There is no singular curvature measure and hence no omitted local-time term.

Conversely, suppose a continuous linear-growth vector and Borel sets satisfy \eqref{eq:EEP}--\eqref{eq:comp3}. The integral equation has locally $W^{1,2}_p$ regularity for finite $p$ and satisfies the PDE with the displayed source. In $\C_i$ its residual is zero; on each $\D_{ij}$ the residual is $-g_{ij}\le0$ and its obstacle is active. It therefore satisfies \eqref{eq:QVI}; uniqueness identifies it with $Q$. The complementarity conditions are essential: value-matching equations alone do not establish uniqueness of proposed boundaries.
\end{proof}

\paragraph{Gaussian evaluation of the integral representation.}
For $u=s-t>0$, write $P_0(u,y;D)=\mathbb P(Y_s^{t,y}\in D)$ and let $P_1$ be the probability under the first-moment-tilted lognormal law. More explicitly, these lognormal laws have log means
\[
 m_0=\log y+(\beta-r-\tfrac12\theta^2)u,\qquad
 m_1=\log y+(\beta-r+\tfrac12\theta^2)u,
\]
and log variance $\theta^2u$. Then
\[
 \E[Y_s^{t,y}\1_D(Y_s^{t,y})]=y e^{(\beta-r)u}P_1(u,y;D).
\]
If $D=(a,b)$, $0\le a<b\le\infty$,
\begin{equation}\label{eq:normal}
 P_\ell(u,y;(a,b))=
 \Phi\left(\frac{\log b-m_\ell}{|\theta|\sqrt u}\right)
 -\Phi\left(\frac{\log a-m_\ell}{|\theta|\sqrt u}\right),\quad\ell=0,1.
\end{equation}
For disjoint interval unions, add these terms; for a general Borel set, integrate the lognormal density over that set. Thus \eqref{eq:EEP} becomes
\begin{align}\label{eq:gaussianEEP}
 Q_i(t,y)={}&\epsilon_i a_r(T-t)y-L_i a_\beta(T-t)\\
 &+\sum_{j\ne i}\left[
  (\epsilon_j-\epsilon_i-rd_{ij})y\int_t^T e^{-r(s-t)}P_1(s-t,y;\D_{ij}(s))\dd s\right.\nonumber\\
 &\hspace{38mm}\left.-(L_j-L_i)\int_t^T e^{-\beta(s-t)}P_0(s-t,y;\D_{ij}(s))\dd s
 \right].\nonumber
\end{align}
Whenever a boundary branch $b_{ij}(t)$ is available, it satisfies
\[
 Q_i(t,b_{ij}(t))=Q_j(t,b_{ij}(t))-d_{ij}b_{ij}(t),
\]
with the functions evaluated through \eqref{eq:gaussianEEP}. There can be several branches associated with a label, or none; the integral formulation does not prescribe their number.

\section{Structure and regularity of the switching boundaries}\label{sec:geometry}
\subsection{Boundary definitions and localisation}\label{sec:pointwise}
We distinguish the boundary of the continuation set from the boundary of a directed contact set. Write
\[
 \Gamma_i=\partial\C_i,\qquad \Gamma_{ij}=\partial\Sset_{ij},
 \qquad
 \Sset_{ij}=\{Q_i=Q_j-d_{ij}y\}.
\]
Boundaries are relative to the time-state domain. The first is the continuation--switching boundary of the current job. The second also records changes between optimal destinations inside a switching region. They need not be the same. Selected sets $\D_{ij}$ additionally depend on a tie-breaking rule, whereas $\Sset_{ij}$ is intrinsic to the value functions.

\begin{proposition}\label{prop:terminal}
For every $0<a<b<\infty$, there is $\eta>0$ such that all jobs are in continuation on $(T-\eta,T]\times[a,b]$. In particular, a sequence of switching points approaching retirement cannot remain in a compact subset of $(0,\infty)$.
\end{proposition}
\begin{proof}
The no-switch policy and the affine reward bound give $|Q_i(t,y)|\le C(T-t)(1+y)$. Write $d_*=\min_{i\ne j}d_{ij}>0$. On $[a,b]$,
\[
 Q_i-Q_j+d_{ij}y\ge d_*a-2C(T-t)(1+b)>0
\]
when $T-t$ is sufficiently small.
\end{proof}

The affine switching gains also give localisation without a regular-point assumption. Let
\[
 g_{ij}(y)=(\epsilon_j-\epsilon_i-rd_{ij})y-(L_j-L_i),\qquad
 K_i=\{y>0:g_{ij}(y)<0\text{ for every }j\ne i\}.
\]

\begin{proposition}\label{prop:pointwise_gain}
If $Q_i(t,y)=Q_j(t,y)-d_{ij}y$ and $(t,y)\in\C_j$, then $g_{ij}(y)\ge0$ at that point, not merely almost everywhere. Consequently,
\begin{equation}\label{eq:core}
 [0,T)\times K_i\subseteq\C_i.
\end{equation}
If all cost triangles are strict, then for every directed pair
\begin{equation}\label{eq:pointwise_localization}
 \Sset_{ij}\subseteq\{(t,y):g_{ij}(y)\ge0\}.
\end{equation}
Thus every actual terminal-destination contact lies in an explicitly known half-line, the whole positive line, or the empty set, according to the signs of the two coefficients of $g_{ij}$.
\end{proposition}
\begin{proof}
Since $j$ is in continuation, $Q_j$ is smooth in a neighbourhood. The function $\phi=Q_j-d_{ij}y$ touches $Q_i$ from below at the contact point. The unconditional viscosity supersolution inequality $-(Q_i)_t-\A Q_i-f_i\ge0$ gives
\[
 0\le-\phi_t-\A\phi-f_i
     =f_j-rd_{ij}y-f_i=g_{ij}(y).
\]
Every contact admits a continuation destination: following active edges cannot form a cycle because all cycle costs are positive, and triangle inequalities collapse the finite chain to a contact with its terminal continuation mode. Therefore a point with all $g_{ij}<0$ cannot be a contact, proving \eqref{eq:core}. Under strict triangle inequalities every active destination is already in continuation, giving \eqref{eq:pointwise_localization} for each pair.
\end{proof}

\subsection{Strict positivity of the switching gain}\label{sec:strictgain}
The positivity condition can be strengthened using the special affine running rewards of this model.

\begin{proposition}\label{prop:strictgain}
Suppose $(t,y)\in\Sset_{ij}\cap\C_j$, with $t<T$ and $y>0$. If $g_{ij}$ is not identically zero, then $g_{ij}(y)>0$. If all triangles are strict for distinct indices, every contact has $g_{ij}>0$; a pair with $g_{ij}\equiv0$ has no contact anywhere before retirement.
\end{proposition}
\begin{proof}
Write $u=Q_i-Q_j+d_{ij}y$ and $g=g_{ij}$. Pointwise nonnegativity was proved in Proposition~\ref{prop:pointwise_gain}. Suppose $g(y_0)=0$ but $g\not\equiv0$. Then $g'$ is a nonzero constant. Waiting in $i$ up to a bounded stopping time $\sigma$ inside a continuation neighbourhood of $j$, then switching to $j$, gives
\begin{equation}\label{eq:strictgain_wait}
 u(t_0,y_0)\ge-\E\int_0^\sigma e^{-\beta s}g(Y_s)\dd s.
\end{equation}
Here time is measured from $(t_0,y_0)$. The inequality follows from dynamic programming, the continuation equality for $Q_j$, and
\[
 \E[e^{-\beta\sigma}d_{ij}Y_\sigma]-d_{ij}y_0
 =-rd_{ij}\E\int_0^\sigma e^{-\beta s}Y_s\dd s.
\]

Let $\xi=\operatorname{sgn}g'$ and stop when $\xi(Y_s-y_0)$ leaves $(-2\delta,\delta)$, or at $M\delta^2$, whichever comes first. Under the scaling $s=\delta^2v$, the stopped process $\xi(Y_{\delta^2v}-y_0)/\delta$ converges to a Brownian diffusion with generator $a_0\partial_{xx}$, where $a_0=\theta^2y_0^2/2$. If $\eta$ is its exit time from $(-2,1)$, solving
\[
 a_0p''(x)=-x,\quad p(-2)=p(1)=0
\]
gives $p(x)=-x^3/(6a_0)+x/(2a_0)-1/(3a_0)$. Hence $\E_0\int_0^\eta X_v\dd v=-1/(3a_0)<0$. Choose a finite $M$ for which the truncated expectation is still negative. Bounded stopped-path convergence then makes the right side of \eqref{eq:strictgain_wait} strictly positive of order $\delta^3$, for small $\delta$. The stopping rule remains in the required neighbourhood and before $T$. This contradicts $u(t_0,y_0)=0$.

It remains to consider $g_{ij}\equiv0$, which means $L_j=L_i$ and $\epsilon_j-\epsilon_i=rd_{ij}$. Follow an optimal policy starting in $j$ up to its first switch, at time $\sigma$, or until $T$ if there is no switch. An alternative starting in $i$ waits until $\sigma$ and then moves directly to the same next job $k$. Relative to switching immediately $i\to j$, its discounted gain is
\[
 e^{-\beta(\sigma-t)}Y_\sigma
 \begin{cases}d_{ij}+d_{jk}-d_{ik},&\sigma<T,\\d_{ij},&\sigma=T.\end{cases}
\]
The preceding wage difference exactly accounts for interest on $d_{ij}$. Under strict triangles all displayed cost coefficients have a positive minimum $\delta_0$, including $k=i$. Therefore
\[
 Q_i(t,y)-Q_j(t,y)+d_{ij}y
 \ge\delta_0\E[e^{-\beta(\sigma-t)}Y_\sigma]
 \ge\delta_0y e^{-r(T-t)}>0.
\]
No contact is possible. Finally, under strict triangles every active destination is in continuation, so the preceding two arguments apply to all contacts.
\end{proof}

\subsection{Scalar reduction and maximal regular components}\label{sec:regularcomponents}
Consider a point $P=(t_0,y_0)\in\partial\C_i$, $0<t_0<T$, $y_0>0$. Suppose one destination $j$ has a strict advantage over all other switching destinations at $P$, and $P\in\C_j$. Both properties persist in a sufficiently small neighbourhood. Put
\begin{equation}\label{eq:u}
 u=Q_i-Q_j+d_{ij}y\ge0.
\end{equation}
There, the entire $N$-mode system reduces for this gap to
\begin{equation}\label{eq:localobstacle}
 \min\{u,-u_t-\A u+g_{ij}(y)\}=0.
\end{equation}
The obstacle is zero, the coefficients are smooth, and $Q_j$ is smooth because it is in continuation.

\begin{proposition}\label{prop:localboundary}
Suppose additionally that $g_{ij}(y_0)>0$ and that $P$ is a regular point of the scalar parabolic obstacle problem: its parabolic blow-up is the nonzero half-space profile
\begin{equation}\label{eq:blowup}
 \frac{u(t_0+\rho^2s,y_0+\rho z)}{\rho^2}
 \longrightarrow\frac{g_{ij}(y_0)}{\theta^2y_0^2}(ez)_+^2
 \quad \text{locally uniformly},
 \qquad e\in\{-1,1\}.
\end{equation}
Here $\rho\downarrow0$, convergence is on compact subsets
of the rescaled $(s,z)$-plane, and the sign $e$ specifies the spatial
side occupied by the continuation region in the limiting profile.
Then the local free boundary is a regular parabolic $C^1$ graph. In one space dimension it can be represented locally as a continuous curve $y=b(t)$. The corresponding wealth boundary $w=\mathcal W_i(t,b(t))$ is continuous.
\end{proposition}
\begin{proof}
Set $s=T-t$, $a=\theta^2/2$, $m=\beta-r-a$, and
\[
 \widehat v(s,\xi)=e^{\beta s}
 u\bigl(T-s,e^{\sqrt a\,\xi-ms}\bigr).
\]
The equation in the positivity set and the vanishing of the relevant
weak derivatives on the zero set give
\[
 \widehat v_s-\widehat v_{\xi\xi}
 =-\widehat g(s,\xi)\1_{\{\widehat v>0\}},\qquad
 \widehat g(s,\xi)=e^{\beta s}
 g_{ij}\bigl(e^{\sqrt a\,\xi-ms}\bigr).
\]
After shrinking the neighbourhood, $\widehat g$ is smooth and
bounded below by a positive constant. Thus this is a standard local
forward parabolic obstacle problem.

Theorem~\ref{thm:regular} and $u(P)=u_y(P)=0$ give the local bound
$0\le u(t,y)\le C(|t-t_0|+|y-y_0|^2)$. The spatial $C^{1,1}$
and temporal derivative estimates also make the rescalings in
\eqref{eq:blowup} locally equibounded and equicontinuous. Their
half-space limit is a separate regular-point hypothesis, not a
consequence of compactness alone.

The changes of variables preserve this spatial half-space type.
The regular-set theorem \cite[Theorem 1.7]{LM} therefore makes the
nearby free boundary a parabolic $C^1$ surface consisting of regular
points. Its nonzero spatial normal yields a local graph $y=b(t)$
in one space dimension. At this stage, parabolic $C^1$ is understood
in the parabolic metric and does not assert ordinary differentiability
of $b$ in time; that stronger conclusion is proved in
Theorem~\ref{thm:cinfty_regular}.
Continuity of the wealth map follows from Theorem~\ref{thm:regular} and \eqref{eq:Rderivatives}.
\end{proof}

We now extend the local graphs along their regular components.
For a pair $i,j$ and $0<t<T$, let
$\mathscr R_{ij}$ consist of those points of $\Gamma_i$ where the active
destination is uniquely $j$, $j$ is in continuation, $g_{ij}>0$, and
the contact gap has the nonzero spatial half-space quadratic blow-up.
By Proposition~\ref{prop:localboundary}, every point of
$\mathscr R_{ij}$ has a neighbourhood in which the free boundary is
a parabolic $C^1$ graph consisting of spatially regular points.
The strict destination inequalities, continuation of $j$, and
positivity of $g_{ij}$ also persist there. Hence $\mathscr R_{ij}$
is relatively open in $\Gamma_i$.

\begin{theorem}\label{thm:regular_components}
Every connected component $G$ of $\mathscr R_{ij}$ is the graph of one
continuous function on a maximal open time interval:
\begin{equation}\label{eq:maximal_graph}
 G=\{(t,b_G(t)):t\in J_G\},\qquad
 J_G\subset(0,T)\text{ an open interval}.
\end{equation}
The graph function $b_G$ is locally parabolic $C^1$ on $J_G$
in the sense of the scalar regular-set theorem. Maximality is relative
to $\mathscr R_{ij}$: the graph cannot be extended as a connected
subset of this regular locus. There are at most
countably many such components. Their wealth images are continuous
graphs $t\mapsto\mathcal W_i(t,b_G(t))$.

A component cannot end at a finite positive state at an interior time
while all the defining regularity conditions persist. More precisely,
every such endpoint cluster point lies outside $\mathscr R_{ij}$.
Thus an endpoint must involve exit from the time-state domain,
loss of uniqueness of the active destination, exit of $j$
from its continuation region, vanishing $g_{ij}$, or failure of the
spatial regular blow-up condition.
\end{theorem}

\begin{proof}
The local graph property makes $\mathscr R_{ij}$ a second-countable
one-dimensional topological manifold, and the projection
$(t,y)\mapsto t$ is a local homeomorphism.
A connected component is homeomorphic either to a circle or to an open
interval. It cannot be a circle, since a circle is compact and
the continuous time projection $(t,y)\mapsto t$ would attain
a maximum $t_*$ at some $(t_*,y_*)\in G$, while the local graph property
gives points of $G$ with time coordinates larger than $t_*$.
Hence $G$ can be
parameterised by an open interval. Along such a
parametrisation, the time coordinate is continuous and locally
injective. It is therefore strictly monotone on the whole parameter
interval: a change of monotonicity would produce a local extremum
and contradict local injectivity. Its image is an open interval $J_G$, and
its inverse gives the single-valued graph \eqref{eq:maximal_graph}.
The local graph representations agree on overlaps, so
$b_G$ inherits the local parabolic $C^1$ regularity throughout $J_G$.

Since $\mathscr R_{ij}$ is a one-dimensional manifold, it is
locally connected, and hence its connected components are open.
Distinct components contain distinct elements of a countable base, so
there are at most countably many. Continuity of the wealth image
follows from continuity of $\mathcal W_i$ and $b_G$.

Finally, a connected component is closed relative to
$\mathscr R_{ij}$. If an interior finite endpoint cluster point belonged
to $\mathscr R_{ij}$, it would therefore belong to $G$. Its time
coordinate would then belong to $J_G$, contradicting that it is an
endpoint of the maximal open interval $J_G$. The listed endpoint
alternatives are exactly the ways to leave $\mathscr R_{ij}$ or the
underlying domain.
\end{proof}

The maximality statement applies to components of the regular locus. Their higher-order smoothness is established next.

\subsection{Smoothness in time and in wealth coordinates}\label{sec:cinfty}
A boundary component is geometrically $C^\infty$ if it is locally an embedded smooth one-dimensional submanifold of the time-state plane. For a graph over time, it suffices that its graph function is $C^\infty$. The following result establishes this property throughout each maximal regular component.

\begin{theorem}\label{thm:cinfty_regular}
Every connected component $G$ of the spatially regular locus $\mathscr R_{ij}$ in Theorem~\ref{thm:regular_components} has
\[
 G=\{(t,b_G(t)):t\in J_G\},\qquad b_G\in C^\infty(J_G).
\]
In particular this is ordinary smoothness in calendar time, and $G$ is a geometrically $C^\infty$ embedded curve. The hypotheses are exactly those defining the regular locus $\mathscr R_{ij}$.
\end{theorem}
\begin{proof}
Fix a point of $G$. In its neighbourhood there is a unique active destination $j$, and $j$ is in continuation. Thus $u=Q_i-Q_j+d_{ij}y$ solves the scalar zero-obstacle problem, with $g=g_{ij}>0$ after shrinking the neighbourhood. Its weak equation is
\begin{equation}\label{eq:scalar_weak_smooth}
 -u_t-\tfrac12\theta^2y^2u_{yy}-(\beta-r)yu_y+\beta u
 =-g(y)\1_{\{u>0\}}.
\end{equation}
On the contact set the relevant weak derivatives vanish almost everywhere. Set
\[
 s=T-t,\quad x=\log y,\quad v(s,x)=u(T-s,e^x),\qquad
 a=\tfrac12\theta^2,\quad m=\beta-r-a.
\]
Then
\[
 v_s-av_{xx}-mv_x+\beta v=-g(e^x)\1_{\{v>0\}}.
\]
The following transformations avoid introducing a time-dependent forcing term. Put
\[
 w(s,x)=e^{mx/(2a)}v(s,x),\qquad c=\beta+\frac{m^2}{4a}>0.
\]
Direct differentiation gives
\begin{equation}\label{eq:gauge_smooth}
 w_s-aw_{xx}+cw=-e^{mx/(2a)}g(e^x)\1_{\{w>0\}}.
\end{equation}
Rescale $x=\sqrt a\,X$ and write $W(s,X)=w(s,\sqrt a\,X)$. Introduce an auxiliary spatial coordinate $z$, restricted to $|z|<\pi/(2\sqrt c)$, and define
\[
 V(s,X,z)=\cos(\sqrt c\,z)W(s,X).
\]
The cosine is positive in this neighbourhood, so
\begin{equation}\label{eq:heat_lift}
 {V_s-V_{XX}-V_{zz}=-F(X,z)\1_{\{V>0\}},\quad V\ge0,}
\end{equation}
where
\[
 F(X,z)=\cos(\sqrt c\,z)e^{mX/(2\sqrt a)}g(e^{\sqrt a X})
\]
is smooth, independent of $s$, and bounded below by a positive constant on a sufficiently small neighbourhood.

The positivity set of $V$ is the product of that of $W$ with the $z$-interval. The previously proved parabolic $C^1$ spatial graph therefore lifts to a parabolic $C^1$ spatial graph in two space dimensions. Corollary 1.4 of Kukuljan \cite{Kuk} states that such a free boundary for \eqref{eq:heat_lift}, with smooth positive spatial forcing, is $C^\infty$ in space and time. The preceding transformations verify these hypotheses for the local equation. Restricting to $z=0$ and undoing the spatial scaling, gauge, and logarithm proves local ordinary $C^\infty$ regularity of $b_G$. The local functions agree on overlaps along the maximal component, which proves the assertion on all of $J_G$.
\end{proof}

\begin{corollary}\label{cor:smooth_traces}
On a regular component, derivatives from the continuation side satisfy
\begin{equation}\label{eq:boundary_curvature}
 u_t=0,\qquad
 u_{yy}=\frac{2g_{ij}(b_G(t))}{\theta^2b_G(t)^2}>0,\qquad
 b_G'(t)=-\frac{u_{ty}}{u_{yy}}(t,b_G(t)).
\end{equation}
If the common utility component $R$ is jointly $C^\infty$ near this component, its wealth image is also a $C^\infty$ curve. This includes the constant-parameter CRRA specification in Subsection~\ref{sec:crra}.
\end{corollary}
\begin{proof}
Smooth boundary Schauder estimates give smooth continuation-side traces of $u$, whose boundary value is zero. Spatial smooth fit gives $u_y=0$. Differentiating $u(t,b_G(t))=0$ gives $u_t=0$, and the continuation equation gives the stated curvature. Differentiating $u_y(t,b_G(t))=0$ gives the velocity formula. These derivatives are taken from the continuation side; no second-order smooth fit across the boundary is being asserted.

At contact the wealth map can be written using the continuing destination:
\[
 \mathcal W_i(t,b_G(t))=-R_y(t,b_G(t))-(Q_j)_y(t,b_G(t))+d_{ij}.
\]
The right side is smooth because $Q_j$ is smooth in its continuation neighbourhood. Its parametrisation has time component $t$, and hence nonzero velocity. The joint smoothness of $R$ is used only for this wealth-coordinate conclusion.
\end{proof}

\subsection{Exact geometry with common labour disutility}\label{sec:globalnew}
When disutility is common across jobs, occupational choice is determined entirely by discounted cash flows. The finite-envelope representation gives a complete geometric classification in this subclass.

\begin{theorem}\label{thm:commonlabor}
Suppose that $L_i=L$ for all $i$, with arbitrary constant wages and a positive effective directed cost matrix. Set
\[
 a=a_r(T-t),\qquad H_i(a)=\max_k\{\epsilon_k a-d_{ik}\}.
\]
Then the exact dual switching value is
\begin{equation}\label{eq:commonlabor}
 {Q_i(t,y)=yH_i(a_r(T-t))-L a_\beta(T-t).}
\end{equation}
Consequently,
\begin{equation}\label{eq:contactcash}
 \Sset_{ij}=\{(t,y):H_i(a_r(T-t))-H_j(a_r(T-t))+d_{ij}=0\}.
\end{equation}
The set of annuity times satisfying this equality is a finite union of closed intervals and isolated points, relative to the available horizon. Its boundary is contained in the combined breakpoints of the two finite affine envelopes $H_i,H_j$.
\end{theorem}
\begin{proof}
The disutility contribution is independent of the policy and equals $-L a_\beta(T-t)$. Dividing the remaining objective by $y$ and changing to the risk-neutral measure gives the expected discounted labour income minus monetary switching costs. By Theorem~\ref{thm:human}, its maximum is $h_i(t)=H_i(a_r(T-t))$, attained by an immediate shortest route to a maximiser followed by no further switching. This proves \eqref{eq:commonlabor} for the stochastic model, not merely for deterministic schedules.

Triangle inequalities imply $H_i\ge H_j-d_{ij}$. On each open cell between consecutive breakpoints of $H_i$ and $H_j$, the nonnegative function $H_i-H_j+d_{ij}$ is affine. It is either identically zero or strictly positive there; an isolated interior zero of a nonnegative affine function is impossible. At a breakpoint it may also vanish. There are only finitely many breakpoints, proving the asserted classification.
\end{proof}

\begin{corollary}\label{cor:common_strict}
If all triangles are strict for distinct indices, then
\begin{equation}\label{eq:common_strict}
 \Sset_{ij}=\left\{(t,y):
  (\epsilon_j-\epsilon_k)a_r(T-t)\ge d_{ij}-d_{ik}
  \ \text{for every }k\right\}.
\end{equation}
Whenever this set is nonempty, its destination $j$ is in continuation. Its annuity-time projection is an interval or a point, and it is therefore connected. Its boundary can nevertheless have two connected components and need not be a graph over calendar time.
\end{corollary}
\begin{proof}
At a contact $i\to j$, if $j$ also contacts $k$, summing the matching equalities would contradict $d_{ik}<d_{ij}+d_{jk}$. Hence $j$ is in continuation and $H_j=\epsilon_j a$. Conversely, if the candidate $\epsilon_j a-d_{ij}$ attains $H_i$, any contact or strictly profitable move from $j$ to $k$ would, by the strict triangle inequality, give a strictly larger candidate for $H_i$. Thus $j$ is in continuation and the candidate gives contact. This proves \eqref{eq:common_strict}. Each displayed condition is a half-line constraint in $a$, so their intersection is an interval, a point, or the empty set. A nondegenerate bounded time interval times $(0,\infty)$ has two boundary components.
\end{proof}

\begin{corollary}\label{cor:common_all_smooth}
If all disutilities are equal, every nonempty connected component of $\Gamma_i$ and $\Gamma_{ij}$ in the interior time-state domain is a straight time line, and is therefore real analytic, for every finite $N$ and every positive effective cost matrix.
\end{corollary}
\begin{proof}
Theorem~\ref{thm:commonlabor} expresses every contact condition through finite piecewise-affine functions of $a_r(T-t)$, independent of $y$. Each contact set, and every finite union of them, is a finite union of time strips and time lines. Its interior-domain boundary is a finite union of disjoint sets $\{t_q\}\times(0,\infty)$. Each connected component has the regular parametrisation $y\mapsto(t_q,y)$. The temporal boundaries in Subsection~\ref{sec:counter3} are smooth in this geometric sense, although they are not spatially regular.
\end{proof}

\begin{proposition}\label{prop:identical}
If $\epsilon_i=\epsilon$ and $L_i=L$ for every job, then
\[
 Q_i(t,y)=\epsilon a_r(T-t)y-La_\beta(T-t),\qquad\Sset_i=\varnothing.
\]
Thus the switching regions may be empty even when every directed move is feasible.
\end{proposition}
\begin{proof}
All policies have the same running reward, and every actual switch strictly reduces the objective by a positive discounted cost. No switching is optimal, and every obstacle inequality is strict.
\end{proof}

\subsection{Disconnected directed contact sets}\label{sec:counter4}
Take
\begin{equation}\label{eq:counter_parameters}
 r=\frac1{10},\quad \beta=\frac15,\quad \theta=\frac3{10},\quad T=5,
 \qquad (\epsilon_1,\epsilon_2,\epsilon_3,\epsilon_4)=(0,1,2,3),
\end{equation}
and $L_1=L_2=L_3=L_4=1$. Use
\begin{equation}\label{eq:counterD4}
 d=\begin{pmatrix}
 0&1&3&6\\
 10&0&5/2&5\\
 10&10&0&4\\
 10&10&10&0
 \end{pmatrix}.
\end{equation}
All off-diagonal entries are positive and all triangle inequalities hold. In particular $d_{14}=d_{12}+d_{24}=6$. These can be the original costs themselves; taking shortest-path closure changes nothing. The four jobs have distinct wages.

Here $0\le a\le a_r(5)=10(1-e^{-1/2})\in(3,4)$. Direct evaluation gives
\begin{align}
 H_1(a)&=\max\{0,a-1,2a-3,3a-6\},\label{eq:H1four}\\
 H_2(a)&=\max\{-10,a,2a-5/2,3a-5\}
         =\max\{a,3a-5\},\label{eq:H2four}\\
 H_3(a)&=2a,\qquad H_4(a)=3a.\nonumber
\end{align}
By \eqref{eq:commonlabor}, these are exact value coefficients for every $y>0$. Therefore the normalised contact gap is
\begin{equation}\label{eq:four_gap}
 \frac{Q_1-Q_2+d_{12}y}{y}=H_1-H_2+1
 =\begin{cases}
 1-a,&0\le a\le1,\\
 0,&1\le a\le2,\\
 a-2,&2\le a\le5/2,\\
 3-a,&5/2\le a\le3,\\
 0,&3\le a\le a_r(5).
 \end{cases}
\end{equation}

\begin{theorem}\label{thm:disconnected_pair}
In the admissible four-job model \eqref{eq:counter_parameters}--\eqref{eq:counterD4},
\begin{equation}\label{eq:S12_four}
 {\Sset_{12}=\{(t,y):a_r(5-t)\in[1,2]\cup[3,a_r(5)],\ y>0\}.}
\end{equation}
This intrinsic pair-contact set has two connected components, and its relative boundary has three connected components. In particular, universal connectedness of $\Sset_{ij}$ is false for an arbitrary positive effective cost matrix.
\end{theorem}
\begin{proof}
The zero set of \eqref{eq:four_gap} is exactly the union in \eqref{eq:S12_four}. Since $a_r(5-t)$ is strictly decreasing and continuous, put
\[
 t_q=5+10\log(1-q/10),\qquad q=1,2,3.
\]
Then $t_3<t_2<t_1$, and
\[
 \Sset_{12}=([0,t_3]\cup[t_2,t_1])\times(0,\infty),\qquad
 \Gamma_{12}=\{t_3,t_2,t_1\}\times(0,\infty).
\]
The intervening open time strip $(t_3,t_2)\times(0,\infty)$ has strictly positive contact gap. Each contact strip is connected, and no connected subset of $\Sset_{12}$ meets both strips because its continuous time projection would have to cross this gap. The same argument separates the three boundary lines.
\end{proof}

Numerically, the exact calendar levels are
\[
 t_3\simeq1.433250561,\qquad t_2\simeq2.768564487,\qquad
 t_1\simeq3.946394843.
\]
On the component with $a\ge3$, job 2 immediately proceeds to job 4. The equal-cost route $1\to2\to4$ and the direct move $1\to4$ give the same payoff. A selector of continuing destinations therefore assigns this component to $\D_{14}$. The intrinsic contact set $\Sset_{12}$ records both optimal routes, whereas the selected set $\D_{12}$ records only moves that continue in job 2.

\subsection{Temporal boundaries under strict triangle costs}\label{sec:counter3}
Retain $r,\beta,\theta,T$ from \eqref{eq:counter_parameters}, and take three jobs with
\begin{equation}\label{eq:counterD3}
 \epsilon=(0,1,2),\qquad L=(1,1,1),\qquad
 d=\begin{pmatrix}0&1&3\\10&0&5/2\\10&10&0\end{pmatrix}.
\end{equation}
Every triangle involving three distinct indices is strict; the smallest slack is $1/2$. The exact coefficients are
\begin{equation}\label{eq:three_h}
 H_1(a)=\max\{0,a-1,2a-3\},\quad
 H_2(a)=\max\{a,2a-5/2\},\quad H_3(a)=2a.
\end{equation}
In this case
\begin{equation}\label{eq:S12_three}
 \Sset_{12}=\{(t,y):1\le a_r(5-t)\le2,\ y>0\}.
\end{equation}
Job 2 is a continuation destination throughout this contact strip, and it is the unique optimal destination in its interior. The strip is connected, while its pair boundary has two disjoint time lines. This example has no equal-cost shortcut degeneracy.

\begin{theorem}\label{thm:temporal_boundary}
Let $t_c=5+10\log(9/10)$, so $a_r(5-t_c)=1$. In a neighbourhood of any $(t_c,y_0)$, $y_0>0$, the actual continuation--switching boundary of job 1 is
\begin{equation}\label{eq:vertical_boundary}
 {\Gamma_1=\{t=t_c\}\times(0,\infty).}
\end{equation}
At each of these points, the active destination is uniquely job 2, job 2 is in continuation, and $g_{12}(y_0)=\tfrac9{10}y_0>0$. Nevertheless, the contact gap $u=Q_1-Q_2+y$ has the temporal blow-up
\begin{equation}\label{eq:temporalblowup}
 {\displaystyle
 \frac{u(t_c+\rho^2s,y_0+\rho z)}{\rho^2}
 \longrightarrow\frac9{10}y_0(s)_+
 }
\end{equation}
locally uniformly. This is not a spatial half-space quadratic blow-up. No neighbourhood of this boundary point admits a representation of the whole boundary as $y=b(t)$. Moreover, $Q_1$ is not differentiable in time at the boundary, although spatial smooth fit holds.
\end{theorem}
\begin{proof}
For $a$ close to 1, $H_1(a)=(a-1)_+$ and $H_2(a)=a$, hence
\begin{equation}\label{eq:three_u}
 u(t,y)=y\,[1-a_r(5-t)]_+.
\end{equation}
For $t>t_c$, the gap is positive and job 1 continues; for $t<t_c$, it switches to job 2. The candidate destination 3 is strictly inferior near $a=1$, since at $a=1$ its entry-adjusted coefficient is $2a-3=-1$, whereas that of job 2 is $a-1=0$. Job 2 does not switch because its next cash threshold is $a=5/2$. The source equals
$g_{12}=(1-rd_{12})y=\tfrac9{10}y$.

Taylor expansion in calendar time gives
\[
 a_r(5-t_c-\rho^2s)=1-\tfrac9{10}\rho^2s+O(\rho^4s^2).
\]
Substitution into \eqref{eq:three_u} proves \eqref{eq:temporalblowup}. The limit is independent of $z$ and is positive for $s>0$, so it cannot equal $c(\nu z)_+^2$ for any $c>0$. At $t=t_c$, the boundary contains every $y$ in a neighbourhood of $y_0$, which excludes a single-valued graph $y=b(t)$.

Finally, differentiating the two explicit formulas on either side yields
\[
 (Q_1)_t(t_c-,y)=e^{-\beta(5-t_c)}-\tfrac9{10}y,
 \qquad (Q_1)_t(t_c+,y)=e^{-\beta(5-t_c)}.
\]
Their difference is nonzero. In contrast, $(Q_1)_y(t_c,y)=0=(Q_2)_y(t_c,y)-1$, exactly the spatial smooth fit relation.
\end{proof}

\begin{figure}[!t]
 \centering
 \includegraphics[width=\linewidth]{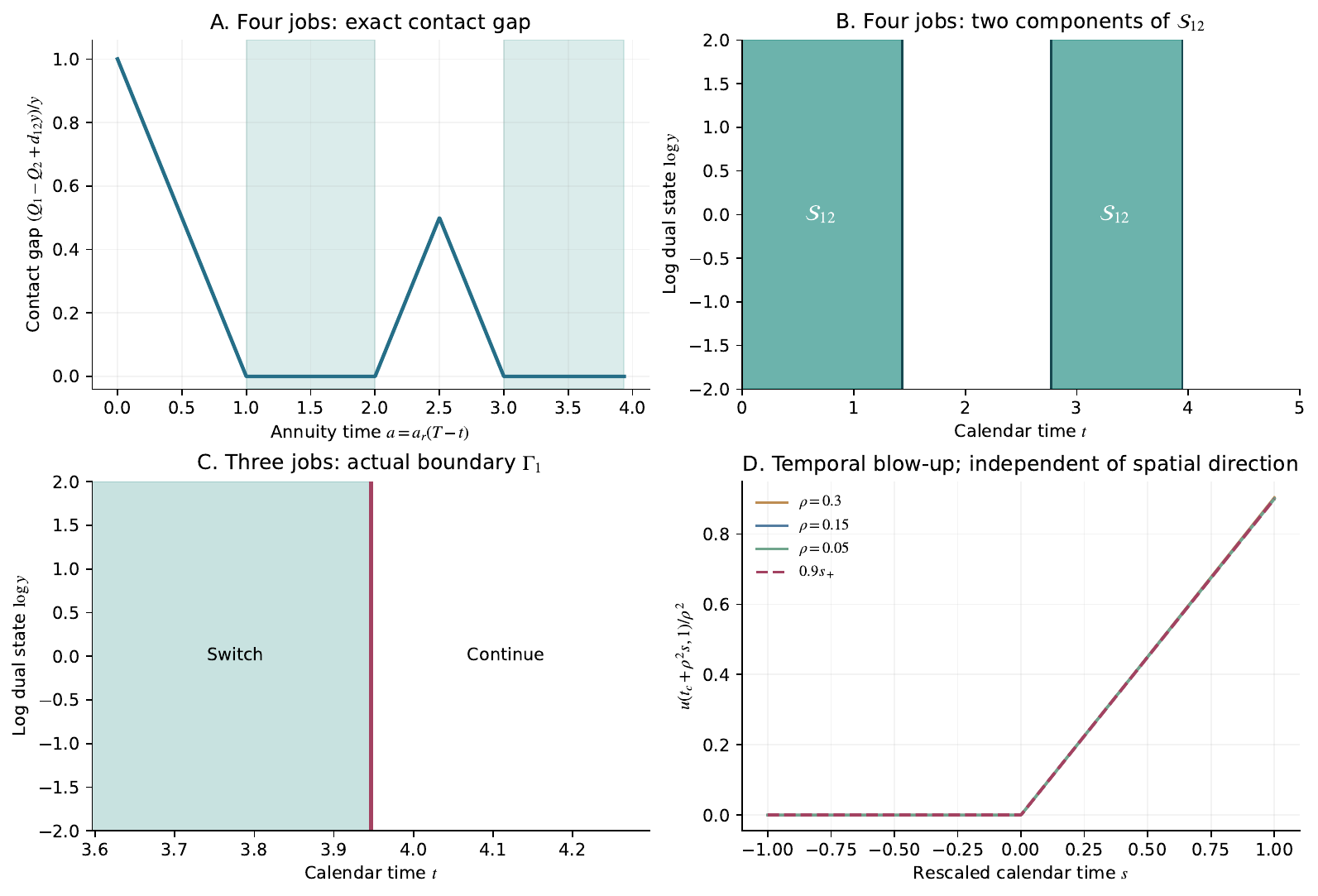}
 \caption{Exact boundary examples. Panels A--B use the four-job model: the normalised contact gap vanishes on two separated annuity-time intervals, giving two components of the intrinsic pair-contact set. Panels C--D use the strict-triangle three-job model: the actual continuation boundary is a time line, and the rescaled contact gap converges to a temporal profile. All curves and regions are evaluated from explicit formulas; no PDE approximation is used.}
 \label{fig:counterexamples}
\end{figure}

The boundary in \eqref{eq:vertical_boundary} is a smooth embedded curve when parameterised by $y$. Its temporal blow-up distinguishes it from the spatially regular graphs of Theorem~\ref{thm:cinfty_regular}, despite strict triangle costs, a unique active destination, and positive switching gain. Figure~\ref{fig:counterexamples} illustrates the exact temporal and disconnected-contact examples.

\paragraph{Financial wealth coordinates.}
Choose admissible CRRA utilities $U_1(c)=-1/c$ and $U_2(w)=-1/w$, so $\gamma=2$, $\kappa=1$. Write
\[
 K=r+\frac{\beta-r}{2}+\frac{\theta^2}{8},\qquad
 A(t)=a_K(T-t)+e^{-K(T-t)}>0.
\]
For either exact common-disutility model,
\begin{equation}\label{eq:explicitwealth}
 w=A(t)y^{-1/2}-H_i(a_r(T-t)),\qquad
 y=\left(\frac{A(t)}{w+H_i(a_r(T-t))}\right)^2.
\end{equation}
These formulas give a time-preserving homeomorphism onto the natural wealth domain. The disconnected set in Theorem~\ref{thm:disconnected_pair} remains disconnected in wealth coordinates. In the three-job example, $H_1(1)=0$, so \eqref{eq:vertical_boundary} becomes $\{t=t_c,w>0\}$. Changing coordinates cannot convert it into a single finite wealth threshold $w=b(t)$.

\subsection{Disconnected switching regions with ordered jobs}
\begin{theorem}\label{thm:disconnected_total}
There is a three-job model with strictly increasing wages, strictly increasing disutilities, and strict triangle costs for which the total switching set $\Sset_2$ and its actual continuation boundary $\Gamma_2$ are both disconnected. One such model uses the market and horizon in \eqref{eq:counter_parameters}, with
\begin{equation}\label{eq:orderedexample}
 \epsilon=(0,1,2),\qquad L=(1,2,5),\qquad d_{ij}=1/10\quad(i\ne j).
\end{equation}
Its time-independent continuation core contains
\begin{equation}\label{eq:orderedcore}
 (0,T)\times\left(\frac{100}{101},\frac{100}{33}\right).
\end{equation}
There are contacts on both sides of this core, so the disconnection concerns the aggregate set $\Sset_2=\Sset_{21}\cup\Sset_{23}$.
\end{theorem}
\begin{proof}
Direct calculation gives
\[
 g_{21}(y)=1-\frac{101}{100}y,\qquad
 g_{23}(y)=\frac{99}{100}y-3.
\]
Both are negative on the interval in \eqref{eq:orderedcore}. Proposition~\ref{prop:pointwise_gain} proves the continuation core and confines contacts $2\to1$ to $y\le100/101$ and contacts $2\to3$ to $y\ge100/33$. These ranges are disjoint, so there are no destination ties in $\Sset_2$.

It remains to prove that both contact sets are nonempty. Let $a=a_r(T)$ and $b=a_\beta(T)>0$. The general value limits are
\[
 \lim_{y\downarrow0}Q_2(0,y)=-b,\qquad
 \lim_{y\uparrow\infty}\frac{Q_2(0,y)}y=2a-1/10,
\]
and $a>1/10$. The exact early-switching-premium representation reads
\begin{align*}
 Q_2(0,y)=ay-2b+\E\int_0^T e^{-\beta s}
 \big[&g_{21}(Y_s)\1_{\Sset_{21}}(s,Y_s)\\
       &+g_{23}(Y_s)\1_{\Sset_{23}}(s,Y_s)\big]\dd s.
\end{align*}
If $\Sset_{21}$ were empty, the remaining nonnegative premium would be at most
$\tfrac{99}{100}\int_0^Te^{-\beta s}\E Y_s\dd s=\tfrac{99}{100}ay$.
Letting $y\downarrow0$ would give $Q_2(0,y)\to-2b$, a contradiction. If $\Sset_{23}$ were empty, the premium would lie between zero and $b$, because $g_{21}\le1$. Dividing by $y$ and letting $y\uparrow\infty$ would give the slope $a$, also a contradiction. Thus there are low-state and high-state contacts.

The vertical line $y=2$ is in continuation at every time and separates these nonempty relatively open parts of $\Sset_2$, so $\Sset_2$ is disconnected. At the time of a low-state contact, move in the spatial slice towards $y=2$ and take the last contact before the continuation core. This gives a point of $\Gamma_2$ below the core. The analogous construction from a high-state contact gives a point above it. Since no boundary point lies inside the open continuation core, the two nonempty parts of $\Gamma_2$ are separated as well.
\end{proof}

\subsection{A corner in an individual boundary component}\label{sec:corner}
The final example combines two separable occupational attributes driven by the same dual state. Each attribute has a smooth boundary near the point of interest, but their intersection creates a corner in the boundary of the joint switching problem.

\paragraph{The leisure attribute.}
Fix $r=1/10$, $\beta=1/5$, $\theta=3/10$, and $T=5$. Attribute $B$ has two states with equal zero wages, disutilities $2$ and $1$, and costs $c^B_{01}=1$, $c^B_{10}=10$. State $1$ dominates state $0$, so there is no optimal reverse switch. Write $q^B_0,q^B_1$ for the two switching values. With remaining horizon $\tau$, define
\begin{equation}\label{eq:B_cost}
 p(\tau,y)=\inf_{0\le\sigma\le\tau}\E\left[
 \int_0^\sigma e^{-\beta s}\dd s
 +e^{-\beta\sigma}Y_\sigma\1_{\{\sigma<\tau\}}\right].
\end{equation}
The choice $\sigma=\tau$ means that no switch occurs. Then
\[
 q^B_1(t,y)=-a_\beta(T-t),\qquad
 q^B_0(t,y)=-a_\beta(T-t)-p(T-t,y),\qquad
 u_B=y-p(T-t,y).
\]

\begin{lemma}\label{lem:B_boundary}
There is a finite positive nonincreasing function $b:(0,T)\to(0,\infty)$ such that
\begin{equation}\label{eq:B_regions}
 \Sset^B_0(t)=\{0<y\le b(t)\},\qquad
 \C^B_0(t)=\{y>b(t)\},\qquad 0<b(t)\le a_\beta(T-t)<5.
\end{equation}
For almost every $t_*$, $b$ is $C^\infty$ on a neighbourhood of $t_*$.
\end{lemma}
\begin{proof}
The objective in \eqref{eq:B_cost} is affine in $y$ for every Brownian stopping rule. Thus $p$ is concave in $y$, and
\[
 0\le p(\tau,y)\le\min\{y,a_\beta(\tau)\}.
\]
The nonnegative function $y-p(\tau,y)$ is convex and vanishes at $y=0$ by continuous extension. Its zero set is consequently an interval starting at zero. Continuity and the displayed bound make its positive endpoint finite whenever a positive contact exists.

The function $p$ is nondecreasing in $\tau$. Indeed, truncate any longer-horizon stopping policy at a shorter horizon and interpret the truncation as no switch. Its integral cost decreases and its switching payment can only be removed. Taking infima proves the claim, and contact sets therefore increase with remaining horizon.

There must be a positive contact for every $\tau>0$. Otherwise monotonicity would give no contact at any shorter horizon either. On that whole horizon the continuation equation and terminal condition would then imply $q^B_0=-2a_\beta(\tau)$ by the linear Feynman--Kac formula, or $p=a_\beta(\tau)$. This contradicts $p\le y$ for small $y$. These observations prove \eqref{eq:B_regions}, including monotonicity in calendar time.

A monotone finite function is differentiable with finite derivative almost everywhere. Fix such a time $t_*$. Then, uniformly for bounded $s$,
\[
 \frac{b(t_*+\rho^2s)-b(t_*)}{\rho}\longrightarrow0.
\]
Thus under parabolic rescaling the contact side contains every compact subset of $\{z<0\}$, whereas the continuation side contains every compact subset of $\{z>0\}$. Here the scalar source satisfies $g_B(b)=1-rb>1/2$. Positive-source nondegeneracy and the one-dimensional blow-up classification \cite{BDM} leave only the positive spatial half-space profile: the other nonzero profiles do not vanish on $\{z<0,s<0\}$. Hence $(t_*,b(t_*))$ is spatially regular. Theorem~\ref{thm:cinfty_regular} now makes the boundary $C^\infty$ near this point. Since the global contact sections are exactly those in \eqref{eq:B_regions}, this is the same function $b$.
\end{proof}

\paragraph{The wage attribute and the four-job model.}
Attribute $A$ has wages $(1,2)$, common disutility $1$, and costs $c^A_{01}=k$, $c^A_{10}=10$, with
\[
 0<k<a_r(T)=10(1-e^{-1/2}).
\]
By the exact common-disutility formula, with $a=a_r(T-t)$,
\begin{equation}\label{eq:A_gap_corner}
 q^A_0=y\max\{a,2a-k\}-a_\beta(T-t),\quad
 q^A_1=2ay-a_\beta(T-t),\quad
 u_A=q^A_0-q^A_1+ky=y(k-a)_+.
\end{equation}
The cash boundary is the time line
\[
 t_*=T+r^{-1}\log(1-rk),\qquad
 \C^A_0=\{t>t_*\},\quad\Sset^A_0=\{t\le t_*\}.
\]
Take the four jobs to be $(0,0),(1,0),(0,1),(1,1)$, labelled $1,2,3,4$. Add the two attributes' wages, disutilities, and costs. Explicitly,
\begin{equation}\label{eq:corner_data}
 \epsilon=(1,2,1,2),\qquad L=(3,3,2,2),\qquad
 d(k)=\begin{pmatrix}
 0&k&1&k+1\\
 10&0&11&1\\
 10&k+10&0&k\\
 20&10&10&0
 \end{pmatrix}.
\end{equation}
Every off-diagonal cost is positive. The triangle inequalities hold because the costs are sums of two directed two-state metrics. A simultaneous change of both attributes costs exactly the sum of its two changes, consistently with the routing convention.

\begin{lemma}\label{lem:product}
The exact four-job value is
\[
 Q_{(a,b)}=q^A_a+q^B_b.
\]
Consequently,
\begin{equation}\label{eq:corner_sets}
 \C_1=\{t>t_*,\ y>b(t)\},\qquad
 \Sset_{14}=\{t\le t_*,\ 0<y\le b(t)\}.
\end{equation}
Job $4$ is in continuation everywhere, and is the unique continuation destination on $\Sset_{14}$.
\end{lemma}
\begin{proof}
Any four-job strategy projects to two attribute strategies. Its running reward and every switching cost split into the sum of their attribute contributions, path by path. Its value is therefore at most $q^A_a+q^B_b$. Conversely, implement optimal attribute strategies together, consolidating simultaneous changes. Costs add exactly, and the combined strategy is admissible and has that value. Probabilistic independence is unnecessary: both decisions use the same Brownian filtration and the same process $Y$.

For job $1$, the gaps to jobs $2,3,4$ are respectively $u_A,u_B,u_A+u_B$. Both attribute gaps are nonnegative. Hence job $1$ continues exactly when both gaps are positive, and it contacts job $4$ exactly when both vanish. Equations \eqref{eq:A_gap_corner} and \eqref{eq:B_regions} prove \eqref{eq:corner_sets}. Job $4$ dominates every other job in flow reward, and its attribute states never switch. On $\Sset_{14}$, jobs $2$ and $3$ still have their respective active attribute moves to job $4$, so neither is in continuation there.
\end{proof}

\begin{figure}[!t]
 \centering
 \includegraphics[width=.8\linewidth]{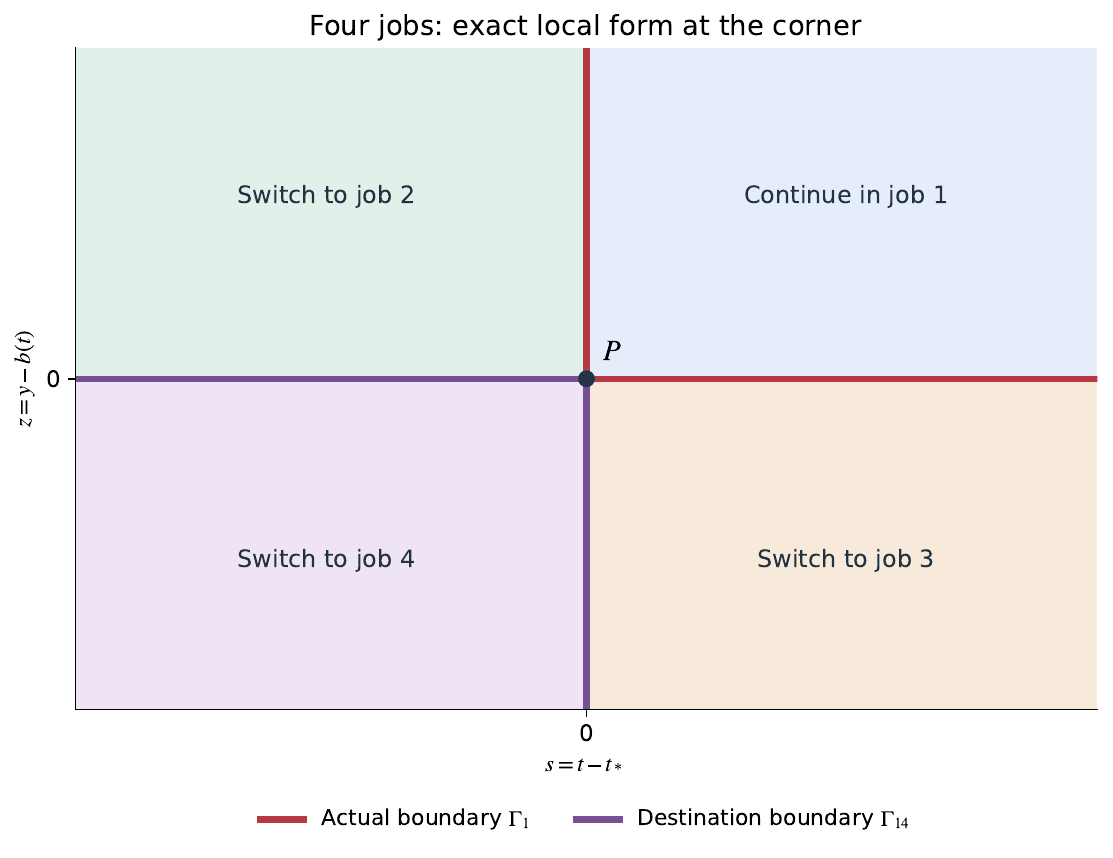}
 \caption{The exact local coordinates \eqref{eq:corner_coordinates}. The red boundary of the current job's continuation region and the purple boundary of the $1\to4$ contact region both have a corner at $P$. The four open quadrants give the indicated terminal destinations. This is the analytic local form \eqref{eq:quadrant_exact}, not a numerical approximation to $b(t)$.}
 \label{fig:geometriccorner}
\end{figure}

\begin{theorem}\label{thm:geometric_corner}
For Lebesgue-almost every $k\in(0,10(1-e^{-1/2}))$, the admissible four-job model \eqref{eq:corner_data} has a connected component of the actual continuation boundary $\Gamma_1$ that is not geometrically $C^1$. A connected component of the individual destination boundary $\Gamma_{14}$ is also not geometrically $C^1$. In particular, both violate geometric $C^\infty$ regularity.

The nonsmooth point is $P=(t_*,b(t_*))$. It has a unique continuation destination, job $4$, and
\[
 g_{14}(b(t_*))=1+\{1-r(k+1)\}b(t_*)>0.
\]
There are also active intermediate destinations at this point.
\end{theorem}
\begin{proof}
By Lemma~\ref{lem:B_boundary}, for almost every $t_*$ the function $b$ is smooth in a neighbourhood of $t_*$. The map $k\mapsto T+r^{-1}\log(1-rk)$ is a smooth diffeomorphism between the specified cost interval and $(0,T)$, with nonvanishing derivative. Thus this condition holds for almost every $k$ in that interval.

Fix such a $k$. Near $P$, use the smooth coordinates
\begin{equation}\label{eq:corner_coordinates}
 s=t-t_*,\qquad z=y-b(t).
\end{equation}
Their Jacobian determinant is one. The exact sets in \eqref{eq:corner_sets} become
\begin{equation}\label{eq:quadrant_exact}
 {\C_1=\{s>0,z>0\},\qquad \Sset_{14}=\{s\le0,z\le0\}.}
\end{equation}
Their local boundaries are respectively the two positive coordinate half-axes and the two negative coordinate half-axes, as shown in Figure~\ref{fig:geometriccorner}. Each is a connected set containing its vertex and two nontrivial incident arcs, so the vertex is not an endpoint of that local component.

Equivalently, in the original coordinates, $\Gamma_1$ contains the arcs
\[
 (t_*,b(t_*)+\eta),\quad\eta\ge0;
 \qquad(t_*+h,b(t_*+h)),\quad h\ge0.
\]
At $P$ their tangent vectors are $(0,1)$ and $(1,b'(t_*))$. Their determinant is $-1$. A $C^1$ embedded curve has a unique tangent line, so this connected component cannot be $C^1$. For $\Gamma_{14}$ the two incident arcs use $\eta\le0$ and $h\le0$ and give the same contradiction. This also shows why a change of smooth coordinates cannot remove the corner.

Finally $k<a_r(5)<4$, so $1-r(k+1)>1/2$. The displayed gain is strictly positive. The assertion about the destination follows from Lemma~\ref{lem:product}. Thus the corner is not caused by vanishing source, zero switching costs, a degenerate diffusion, or a boundary at $y=0$ or $T$.
\end{proof}

The two incident arcs belong to the same connected component. The example therefore gives a local failure of geometric regularity in both $\Gamma_1$ and $\Gamma_{14}$. Its additive costs include triangle equalities; the role of strict triangle costs in excluding such multiple-destination corners remains a separate question.

\begin{corollary}\label{cor:wealth_corner}
For the CRRA utilities $U_1(c)=-1/c$ and $U_2(w)=-1/w$, the four-job example also has nonsmooth individual boundary components in financial wealth coordinates.
\end{corollary}
\begin{proof}
Let $A(t)>0$ be the smooth CRRA coefficient in \eqref{eq:explicitwealth}. On $\Sset_{14}$ spatial smooth fit and $Q_4=2a_r(T-t)y-2a_\beta(T-t)$ give
\[
 \mathcal W_1(t,y)=A(t)y^{-1/2}-2a_r(T-t)+k+1.
\]
The right side defines a smooth ambient map $(t,y)\mapsto(t,w)$ with strictly negative $y$-derivative. It is a local diffeomorphism and therefore preserves the corner of $\Gamma_{14}$.

For the actual boundary $\Gamma_1$, the time-line arc remains a time-line arc, with nonzero spatial velocity. On the other arc, where $t>t_*$ and $y=b(t)$, attribute $A$ continues in state $0$ and attribute $B$ has smooth fit to state $1$, giving
\[
 w=A(t)b(t)^{-1/2}-a_r(T-t)+1.
\]
This is a smooth graph with finite velocity and time component one. Its tangent at the vertex cannot be parallel to the time-line arc. Hence the actual boundary has a wealth-coordinate corner as well.
\end{proof}

\section{Implementation and economic implications}\label{sec:computation}
\subsection{CRRA preferences and optimal controls}\label{sec:crra}
Take $\gamma>0$, $\gamma\ne1$, and $\kappa>0$, with
\[
 U_1(c)=\frac{c^{1-\gamma}}{1-\gamma},\qquad
 U_2(T,w)=\kappa\frac{w^{1-\gamma}}{1-\gamma}.
\]
Set
\begin{equation}\label{eq:K}
 K=r+\frac{\beta-r}{\gamma}+\frac{\gamma-1}{2\gamma^2}\theta^2,
 \qquad A(t)=a_K(T-t)+\kappa^{1/\gamma}e^{-K(T-t)}.
\end{equation}
The formula is valid for any sign of $K$, using $a_0(u)=u$. Direct lognormal integration gives
\begin{equation}\label{eq:CRRA}
 R(t,y)=A(t)\frac{\gamma}{1-\gamma}y^{(\gamma-1)/\gamma},
 \quad -R_y=A(t)y^{-1/\gamma},
 \quad R_{yy}=\frac{A(t)}\gamma y^{-1/\gamma-1}.
\end{equation}
Therefore the single scalar equation
\begin{equation}\label{eq:CRRAinvert}
 w=A(t)y^{-1/\gamma}-(Q_i)_y(t,y)
\end{equation}
determines the unique multiplier, after which
\begin{equation}\label{eq:CRRAcontrols}
 c^*=y^{-1/\gamma},\qquad
 \pi^*=\frac{\theta}{\sigma}\left\{\frac{A(t)}\gamma y^{-1/\gamma}
                                      +y(Q_i)_{yy}(t,y)\right\}.
\end{equation}
For a post-retirement Merton consumption horizon $D$, the choice $\kappa^{1/\gamma}=a_K(D)$ gives $A(t)=a_K(T+D-t)$ and recovers the retirement continuation specification of Yang and Jeon \cite{YJ}. In the logarithmic case, the same control formulas hold with $\gamma=1$ and $A(t)=a_\beta(T-t)+\kappa e^{-\beta(T-t)}$.

The formulas yield several economic implications. First, adding available jobs weakly increases the dual switching value and maximal human capital, because every old strategy remains available. It also weakly increases primal utility at a fixed feasible wealth. Second, reducing any original switching cost weakly improves the attainable choice set in monetary terms. Third, the dual job regions depend on $(r,\beta,\theta,\epsilon,L,d,T)$, while risk aversion changes their wealth images through \eqref{eq:CRRAinvert}. Fourth, when $\theta>0$, the switching curvature term in \eqref{eq:CRRAcontrols} is nonnegative at a fixed dual state. This last comparison holds at the same $y$; it does not by itself order portfolios at the same financial wealth, because the multiplier also changes.

\subsection{Time discretisation and convergence}\label{sec:discretisation}
For an integer $M\ge1$, let $\Delta=T/M$ and $t_m=m\Delta$, $m=0,\ldots,M$, and set $q_i^M(y)=0$. Let $\mathsf P_\Delta$ denote the exact, undiscounted transition semigroup of $Y$. Work backwards using
\begin{align}
 B_j^m(y)&=\epsilon_j a_r(\Delta)y-L_j a_\beta(\Delta)
           +e^{-\beta\Delta}(\mathsf P_\Delta q_j^{m+1})(y),\label{eq:Bgrid}\\
 q_i^m(y)&=\max_{j\in\I}\{B_j^m(y)-d_{ij}y\}.\label{eq:grid}
\end{align}
The maximising $j$ is the selected job for the next interval. Including $j=i$ allows continuation. The matrix $d$ incorporates every instantaneous multi-step route. The cost per step is $O(N^2)$ comparisons per spatial node, and the continuous state dimension stays one.

\begin{theorem}\label{thm:grid}
Before spatial approximation, \eqref{eq:Bgrid}--\eqref{eq:grid} is the exact value of the problem in which switches are restricted to the grid. On nested grids the values increase to $Q$ as $\Delta\downarrow0$, locally uniformly in the dual state at common dates. Interpolation in the horizon extends convergence locally in time. No preassigned switching topology is needed.
\end{theorem}
\begin{proof}
Backward conditioning gives the recursion. A grid policy is an admissible continuous-time policy, so its value is at most $Q$. Nested grids enlarge the allowed policy set. Given an $\varepsilon$-optimal continuous-time policy, first truncate it at a deterministic number $m$ of switches, making its payoff loss arbitrarily small by integrability. Round each of these at most $m$ stopping times up to the next grid date. Multiple switches assigned to one date can be replaced by a shortest route; this weakly reduces their monetary cost. A route rounded to $T$ can be deleted. Continuity of $Y$, bounded wage and disutility rates, and the integrable bound $m d_{\max}\sup_{s\in[t,T]}Y_s$, where $d_{\max}=\max_{i,j}d_{ij}$, give convergence of the rounded payoff to at least the original truncated payoff. Thus the limiting grid value is at least $Q-\varepsilon$. Let $\varepsilon\downarrow0$. Convexity gives compact equicontinuity in $y$, and the horizon truncation estimate on common grid dates supplies time equicontinuity. This proves the stated local convergence.
\end{proof}

The numerical implementation accepts arbitrary finite arrays of wages and disutilities and a positive directed cost matrix. It computes directed shortest paths and uses \eqref{eq:grid}. The lognormal expectation of a piecewise linear interpolant in $y$ is integrated analytically with normal distribution functions. The only spatial approximation is the interpolation of the value function and its extension beyond the computational interval. The finite spatial endpoints are extended linearly, anchored at the computed endpoint value, with the cash-maximising slope at large $y$ and the minimum-disutility cash slope at small $y$. Domain expansion is checked separately.

\subsection{Five-job illustration}
The parameters, in illustrative units, are
\[
 T=30,\quad r=0.02,\quad\beta=0.04,\quad\theta=0.30,
\]
\[
 (\epsilon_1,\ldots,\epsilon_5)=(0.3,0.7,1.1,1.5,0.9),\qquad
 (L_1,\ldots,L_5)=(0.3,0.7,1.4,2.7,1.8),
\]
with original and effective costs
\begin{equation}\label{eq:numericcost}
 k=\begin{pmatrix}
 0&1.8&4.8&6.5&0.6\\0.8&0&1.3&3.8&0.9\\1&0.7&0&2.4&1.6\\
 1.2&1&0.8&0&1.4\\1.1&0.5&1.4&2&0
 \end{pmatrix},\quad
 d=\begin{pmatrix}
 0&1.1&2&2.6&0.6\\0.8&0&1.3&2.9&0.9\\1&0.7&0&2.4&1.6\\
 1.2&1&0.8&0&1.4\\1.1&0.5&1.4&2&0
 \end{pmatrix}.
\end{equation}
For example, the route $1\to5\to4$ costs $2.6$, compared with the original direct charge $6.5$. Job 5 has a lower wage and greater disutility than job 3, but its low entry cost and route position distinguish it in the switching problem. A wage-disutility comparison alone is insufficient for deleting its state.

\begin{figure}[p]
 \centering\includegraphics[width=.92\linewidth,height=.82\textheight,keepaspectratio]{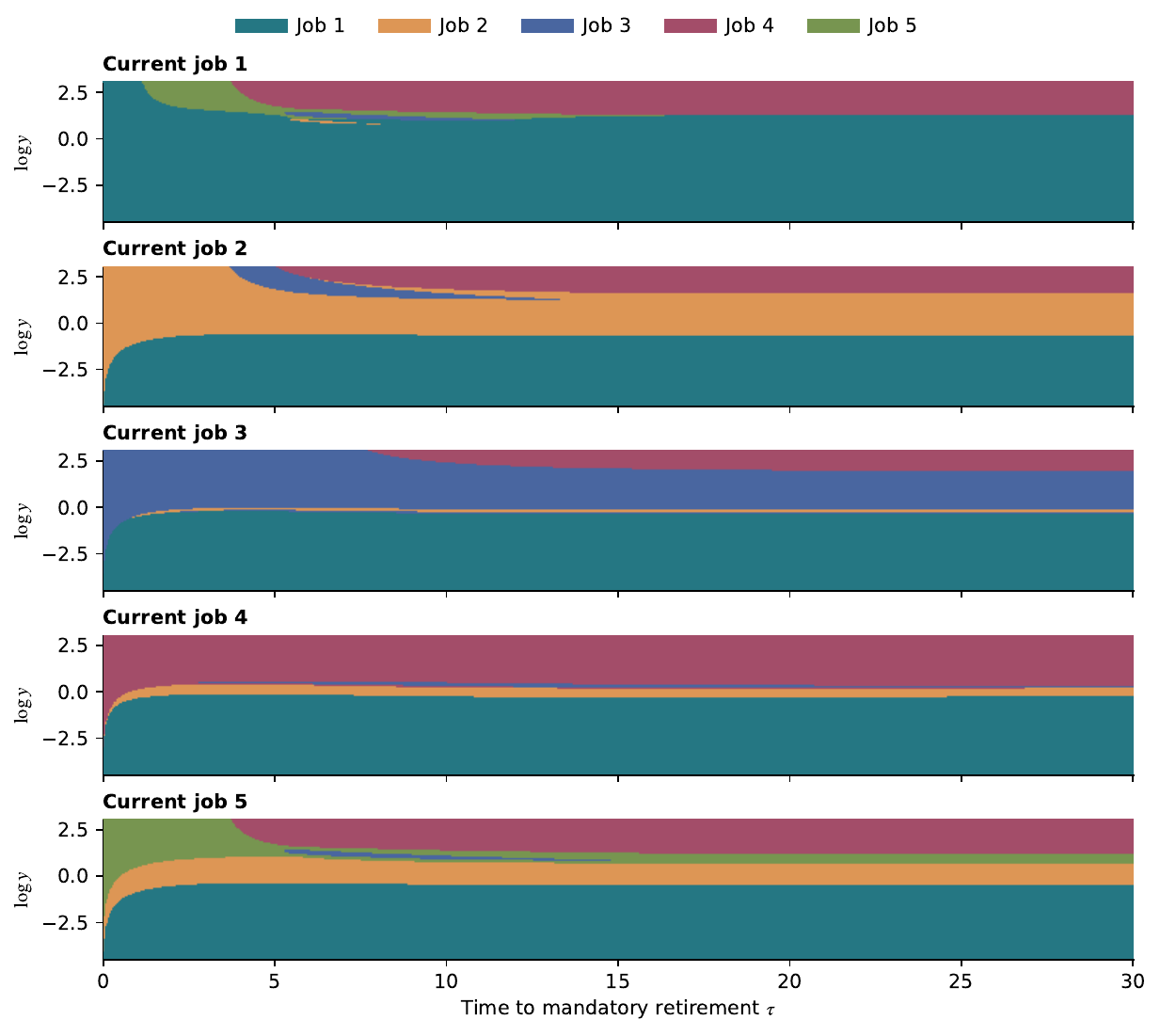}
 \caption{Numerical grid decisions for five arbitrary jobs. Each panel fixes the current job; the colour identifies the job selected for the next interval. The horizontal axis is remaining time, increasing to the right. Selecting the current job means continuation. Effective shortest routes are charged. The plot displays only the interior of the larger computational state interval. The decisions are obtained from the time-grid switching recursion.}
 \label{fig:policy}
\end{figure}
\begin{figure}[htbp]
 \centering\includegraphics[width=\linewidth]{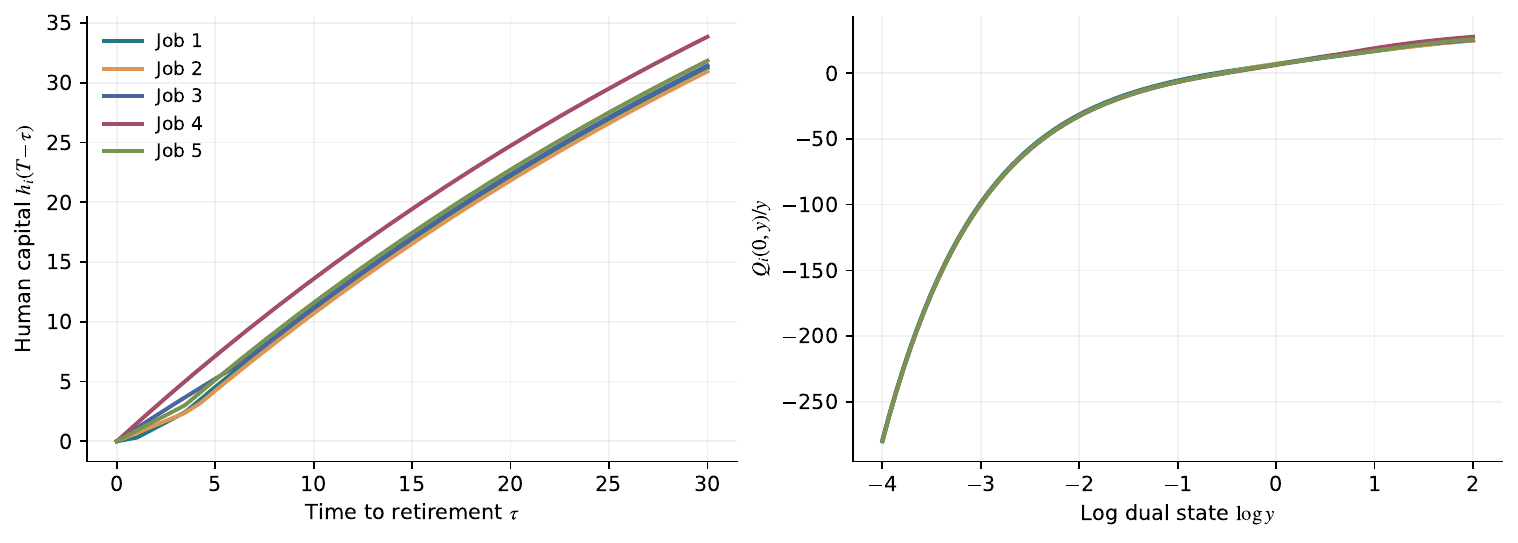}
 \caption{Exact human-capital envelopes (left) and computed normalised dual switching values at initial time (right). Different initial jobs retain different options and switching liabilities.}
 \label{fig:values}
\end{figure}

Figure~\ref{fig:policy} displays the resulting occupational decisions, and Figure~\ref{fig:values} compares the human-capital envelopes and dual values. The main grid has 1,200 time intervals and 1,201 log-state nodes on $[-8,6]$. The coarse and fine grids use $(600,601)$ and $(2400,2401)$ time intervals and state nodes, respectively. Table~\ref{tab:numerics} reports the identical-job benchmark error, the intervention gap, and interior convexity. It also compares values at $y=0.05,0.2,1,3$ across these grids and after expanding the log-state domain to $[-10,8]$. These are value-level diagnostics; they do not establish convergence rates for individual switching curves.
\begin{table}[htbp]
\centering
\caption{Value-function diagnostics for the five-job computation.}\label{tab:numerics}
\small
\begin{tabular}{lr}
\toprule
Diagnostic & Observed value\\
\midrule
Identical-jobs maximum absolute error & $1.67\times10^{-10}$\\
Minimum intervention gap & $-1.59\times10^{-12}$\\
Minimum interior slope increment in $y$ & $6.08\times10^{-3}$\\
Coarse-to-medium maximum probe difference & $5.91\times10^{-2}$\\
Medium-to-fine maximum probe difference & $2.97\times10^{-2}$\\
Expanded-domain maximum probe difference & $3.56\times10^{-4}$\\
\bottomrule
\end{tabular}
\end{table}
The intervention violation is at floating-point roundoff scale. Refinement reduces the observed value differences; it does not furnish a certified error bound for the continuous-time solution or its switching boundaries.

\subsection{Vanishing switching costs}
\begin{proposition}\label{prop:zerocost}
Let $d_{ij}^{\eta}=\eta d_{ij}$ with $\eta>0$. As $\eta\downarrow0$, the switching value increases to an initial-job-independent limit
\begin{equation}\label{eq:zerocost}
 Q^0(t,y)=\E\int_t^T e^{-\beta(s-t)}
                     \max_{j\in\I}(\epsilon_jY_s^{t,y}-L_j)\dd s.
\end{equation}
On compact initial-state sets, an elementary bound is
\[
 0\le Q^0(t,y)-Q_i^\eta(t,y)\le C(1+y)\eta^{1/3},\qquad 0<\eta\le1,
\]
with a constant depending on the fixed parameters and horizon.
\end{proposition}
\begin{proof}
The pointwise maximum bounds every running reward and removing positive costs gives the upper bound. On a grid of width $\Delta$, choose a running-reward-maximising job at each left endpoint and hold it until the next date. The functions $f_j$ are uniformly Lipschitz in $y$, and geometric Brownian increments have expected size at most $Cy\sqrt\Delta$ over a grid interval. Integrating the resulting reward loss over the horizon gives $C(1+y)\sqrt\Delta$. At most $1+(T-t)/\Delta$ switches are made, and their expected discounted cost is at most $C(1+y)\eta/\Delta$. Choosing $\Delta$ of order $\eta^{2/3}$ proves the bound; the last partial interval and an initial switch are absorbed into the constant. Monotonicity follows from decreasing costs.
\end{proof}

At zero cost, tracking an instantaneous reward maximiser may require infinitely many switches near Brownian crossings of indifference levels. Formula~\eqref{eq:zerocost} thus gives the limit of the positive-cost values; attainment at zero cost depends on the admissible class.

\subsection{The two-job benchmark}
For $N=2$, relabel the jobs as $0,1$ and set $d_{01}=\zeta_0$, $d_{10}=\zeta_1$. Shortest-path closure leaves these costs unchanged. Under $\epsilon_0<\epsilon_1$, $L_0<L_1$,
\[
 h_1(t)=\epsilon_1a_r(T-t),\qquad
 h_0(t)=\max\{\epsilon_0a_r(T-t),\epsilon_1a_r(T-t)-\zeta_0\}.
\]
If $r\zeta_0<\epsilon_1-\epsilon_0$, the two candidates in the maximum defining $h_0$ coincide at the critical remaining time
\[
 T_1=-\frac1r\log\left(1-\frac{r\zeta_0}{\epsilon_1-\epsilon_0}\right).
\]
Thus immediate switching increases job 0's net human capital exactly
when $T-t>T_1$. At $T-t=T_1$, it is the values of staying in job 0
and paying to enter job 1 that agree; $h_1-h_0=\zeta_0>0$ there.
The two QVIs in \eqref{eq:QVI} reduce to the scalar double-obstacle problem for $P=Q_1-Q_0$ with obstacles $-\zeta_1y$ and $\zeta_0y$. Under the two-job boundary theorem of \cite{YJ}, the selected sections become a lower half-line for $1\to0$ and an upper half-line for $0\to1$. Substitution into \eqref{eq:gaussianEEP} recovers their two coupled integral equations, and \eqref{eq:CRRAcontrols} recovers their optimal-control formulas.

For general $N$, the corresponding closed object is the interconnected system with contact-set complementarity. The number of jobs changes the discrete comparison problem, while the complete-market dual state remains one-dimensional.

\section{Conclusion}\label{sec:conclusion}
Occupational flexibility before a fixed retirement date can be incorporated into consumption and investment through a one-dimensional dual state and finitely many interconnected value functions. Directed shortest routes determine the exact natural borrowing limit. The affine occupational rewards then yield spatial $C^{1,1}$ regularity, strong duality, and explicit recovery of consumption and portfolio policies without imposing a threshold structure. The switching QVI, its obliquely reflected BSDE representation, and the set-based integral formula provide complementary descriptions of the same value vector.

The boundary results distinguish value regularity from geometric regularity. Directed contact sets may be disconnected, and multiple active destinations can create a corner in an individual boundary component. The maximal $C^\infty$ graph theorem therefore concerns the spatially regular locus, whereas spatial smooth fit and control recovery apply to all contacts. The time-grid scheme retains the full destination structure, and the exact examples provide benchmarks for assessing numerical switching decisions.

Further questions include geometric classifications under strict triangle costs, occupational models with unhedgeable income risk, and learning-based approximations when the model is not known. The reinforcement-learning formulation of Huang et al.\ \cite{HLYZ} suggests one possible direction, but adapting its analysis to unbounded affine rewards and triangle equalities requires additional work. Neither the value convergence proved here nor entropy regularisation alone establishes convergence of individual switching boundaries.

\end{document}